\documentclass[preprint,11pt]{elsarticle}     
\usepackage{mathptmx}      
\usepackage{latexsym}
\usepackage{graphicx,graphics}
\usepackage{amsmath,amssymb,latexsym,epsfig} 
\usepackage[english]{babel}
\usepackage{amsfonts,amsmath}
\usepackage{graphicx,psfrag}
\usepackage{booktabs,caption}
\usepackage{multirow}
\usepackage{subfigure}
\usepackage{amsthm}
\usepackage{epstopdf}
\usepackage{algorithm}
\usepackage{algorithmic}
\usepackage{xcolor}
\usepackage{listings}

\def \x {{\bf x}}
\def \u {{\bf u}}
\def \v {{\bf v}}
\def \k {{\bf k}}
\def \w {{\bf w}}
\def \bu {{\mathbf u}}

\def\RR{\mathbb{R}}

\def\CC{\mathbb{C}}
\def \pmatrix{ \left( \begin{array} }
\def \endpmatrix{ \end{array} \right) }

\definecolor{mygreen}{RGB}{40, 155, 120}

\newcommand{\eeq}{\end{equation}}
\newcommand{\beq}{\begin{equation}}

\newtheorem{remark}{Remark}[subsection]
\begin{document}

\title{Matrix-oriented discretization methods for reaction-diffusion PDEs: comparisons and applications}

%

\begin{frontmatter}

\author{Maria Chiara D'Autilia, Ivonne Sgura\corref{cor1}}
\ead{mariachiara.dautilia@unisalento.it, ivonne.sgura@unisalento.it}
\address{Dipartimento di Matematica e Fisica ``E. De Giorgi'',
Universit\`{a} del Salento, Via per Arnesano, 73100 Lecce,
Italy}
\cortext[cor1]{Corresponding author}

\author{Valeria Simoncini}
\ead{valeria.simoncini@unibo.it}
\address{Dipartimento di Matematica, Alma Mater Studiorum - Universit\`{a} di Bologna, Piazza di Porta San Donato  5, I-40127 Bologna, Italy, and IMATI-CNR, Pavia.}


%
%
%

\begin{abstract}
Systems of reaction-diffusion partial differential equations (RD-PDEs) are widely applied for 
modelling life science and physico-chemical phenomena. In particular, the coupling between
diffusion and nonlinear kinetics can lead to the so-called Turing instability, giving
rise to a  variety of spatial patterns (like labyrinths, spots, stripes, etc.) attained 
as steady state solutions for large time intervals. To capture the morphological peculiarities of 
the pattern itself, a very fine space discretization may be required, 
limiting the use of standard (vector-based) ODE solvers in time because of excessive computational costs.
{We show that the structure of the diffusion matrix can be exploited so as to use matrix-based 
versions of time
integrators, such as Implicit-Explicit (IMEX) and exponential schemes. This implementation entails
the solution of a sequence of discrete matrix problems of significantly smaller dimensions than in the 
vector case, thus allowing for a much finer problem discretization. 
We illustrate our findings by numerically solving  the Schnackenberg model, prototype of 
RD-PDE systems with Turing pattern solutions, and the DIB-morphochemical model describing 
metal growth during battery charging processes.}\end{abstract}

\begin{keyword}
 Reaction-diffusion PDEs \sep Turing patterns \sep IMEX methods \sep ADI method \sep Sylvester equations \sep Schnackenberg model
\end{keyword}

\end{frontmatter}

\section{Introduction}
We are interested in the numerical solution of reaction-diffusion
partial differential equations (RD-PDEs) of the type
\begin{eqnarray}\label{eqn:main1}
u_t = \ell(u) + f(u), \quad u=u(x,y,t), \quad with \quad (x,y)\in
\Omega \subset \mathbb{R}^2,\quad t \in ]0, T],
\end{eqnarray}
{with given initial condition $u(x,y,0)=u_0(x,y)$} and
appropriate boundary conditions on the spatial domain. We assume
that the diffusion operator $\ell$ is linear in $u$, while the
function $f$ contains the nonlinear {reaction} terms.
This setting can be generalized to the system case, which reads as
\begin{equation}\begin{cases}\label{eqn:main2}
u_t= \ell_1(u) + f_1(u,v),  \\
v_t= \ell_2(v) + f_2(u,v), \quad with \quad (x,y)\in \Omega
\subset \mathbb{R}^2,\quad t \in ]0, T]. \end{cases}
\end{equation}
RD-PDE systems describe mathematical models of interest in many classical scientific fields 
like chemistry \cite{Vanag, wit}, biology \cite{Murray, Maini}, ecology \cite{malchow, Sherratt}, 
but also in recent applications concerning for example metal growth by electrodeposition 
\cite{TuringJCAM2012, EJAM2015, Lacit2017}, 
tumor growth \cite{Chaplain}, biomedicine \cite{Gerisch} and cell motility \cite{Madzva1}.
In particular, since the pioneering work of Alan Turing at the origin of mathematical description of morphogenesis \cite{Turing52},
it has been shown that the coupling between diffusion and nonlinear kinetics can lead to the so-called \emph{diffusion-driven}
 or \emph{Turing} instability, giving rise to a wide variety of spatial patterns (like labyrinths, spots, stripes, etc.)
as stationary solutions to \eqref{eqn:main2}; see, e.g.,
\cite{Murray,been, boz, EJAM2015}. 
These models have been
classically solved in a planar domain $\Omega$, however more
recently certain physical applications have motivated the solution
of \eqref{eqn:main2} on stationary or time evolving surfaces where
the diffusion $\ell(u)$ is defined in terms of the
Laplace-Beltrami operator \cite{Frittelli1,Frittelli2,Madzva2}.
The numerical treatment of these models require a space
discretization, and a time integration, commonly performed by
Implicit-Explicit schemes \cite{Ascher95, Ruuth95}, or the
Alternating Direction Integration (ADI) approaches \cite{EJAM2015,
TuringJCAM2012,JCAM2016}. {It is well known that } the Method of
Lines (MOL) {based on classical semi-discretizations in space
(e.g. finite differences, finite elements)
 rewrites (\ref{eqn:main1}) as an ODE system}
\begin{eqnarray} \label{eqn:main1_discr}
\dot \u = A \u + f(\u), \quad \u(0) = \u_0
\end{eqnarray}
where the entries of the matrix $A$ stem from the discretization
of the spatial derivatives involved in the diffusion operator $
\ell(u)$, while the vector $\bu$ contains the coefficients for the
approximation of the sought after function $u$, in the chosen
basis. Analogous ODE equations are obtained by the semi-discretization of
the PDE system \eqref{eqn:main2}, that is
\begin{equation}\label{eqn:main2_discr} \begin{cases}
\dot \u = A_1 \u + f_1(\u,\v) , \quad \u(0) = \bu_0, \\
\dot \v = A_2 \v + f_2(\u,\v),\quad \v(0) = \v_0.\end{cases}
\end{equation}

We are interested in exploring time stepping strategies that can efficiently handle the nonlinear part, by judiciously exploiting
the coefficient matrix structure in the linear part.
In particular, we discuss the situation where the given domain is sufficiently regular so that the linear
differential operator $\ell(u)$ can be discretized by means of a
tensor basis; for instance, this is the case for finite difference
methods, or for certain finite element techniques or spectral
methods. In this framework, the physical space can be mapped into a so-called ``logical space'', typically represented by a
rectangle; see, e.g., \cite{Knupp.Steinberg.92}. To simplify the
presentation, and to adhere to the application we are going to
focus on, in the following we shall restrict the discussion to a
rectangular domain, say $\Omega = [0, \ell_x] \times [0,
\ell_y]$. With these premises, the discretization of the linear
diffusion operator leads to a matrix $A$ of the form
\begin{eqnarray}\label{eqn:kron}
A = I \otimes T_1 + T_2^T \otimes I \in\RR^{N_xN_y\times N_xN_y},
\end{eqnarray}
where $\otimes$ is the Kronecker operator, and $T_1$ ($T_2$) contains
the approximation of the second order 
derivative\footnote{In practice, the two matrices
also contain information on the PDE boundary conditions; See for 
example \cite{TuringJCAM2012} {and Section 4 in the case of zero Neumann BCs.}}
in the
$x$-direction ($y$-direction):  $N_x$ and $N_y$ are the numbers of mesh interior nodes
in the x- and y-directions, respectively.
For example, in the case of finite differences $h_x = \ell_x/(N_x+1), \ h_y = \ell_y/(N_y+1)$ will be the
corresponding space meshsizes.

With this hypotheses, at each time $t\in [0, T]$ it is possible to explicitly employ the
matrix $U(t) \in \RR^{N_x\times N_y}$ containing the same
components of $\bu(t)$, with $U_{i,j}(t) \approx u(x_i,y_j,t)$,
that is, the rows and columns of $U$ explicitly reflect the space
grid discretization of the given problem. 
{The vector} $\bu$ corresponds to the {\tt vec} operation of the matrix $U$, where
each column of $U$ is stuck one after the other. In a finite difference discretization
this {implements} a lexicographic order of the nodes in the rectangular grid.
With this notation, for $A$ in (\ref{eqn:kron}) we {have} 
that $A \bu = {\rm vec}(T_1 U + U T_2)$.
Then (\ref{eqn:main1_discr}) can be written as the following differential matrix equation
\begin{eqnarray}\label{eqn:matrixform}
\dot U = T_1 U + U T_2 + F(U), \quad U(0)=U_0,
\end{eqnarray}
where $F$ is the nonlinear vector function $f(\bu)$ evaluated
componentwise, and ${\rm vec}(U_0)=\bu_0$ is the initial
condition.
Analogously, under the same discretization framework, the system
in \eqref{eqn:main2_discr} can be brought to the matrix form
\begin{equation}\label{eqn:matrixform_sys} \begin{cases}
\dot U = T_{11} U + U T_{12} + F_1(U,V), \quad U(0)=U_0, \\
\dot V = T_{21} V + V T_{22} + F_2(U,V), \quad V(0)=V_0 \end{cases}
\end{equation}
with obvious notation for the introduced quantities.

{These matrix forms provide a quite different perspective at
the time discretization level than classical approaches, allowing
to significantly reduce the memory and computational requirements.
{The use of matrix-based approaches has only very recently been explored
in a systematic manner, in the context of linear or quadratic
matrix terms, such as the Sylvester and Riccati differential equations,
see, e.g., \cite{Benner2018,Beher2018,Stillfjord2015}. 
In particular, matrix ODE equations like
\eqref{eqn:matrixform} -
in the case when the term $F(U)$ allows for low-rank approximations -
have been considered
in \cite{Mena2018}; in that article,} the authors propose a low-rank strategy
for matrix approximation in time based on Lie-Trotter and Strang
splittings.}
Here, to complete the PDE approximation in the two cases above
\eqref{eqn:matrixform} and \eqref{eqn:matrixform_sys}, we show
that the standard time discretization strategies can be tailored
to the matrix equation setting, with several numerical advantages.

{The paper is structured as follows.} After a brief survey in section~\ref{sec:classical} of methods
usually applied to discretize in time (\ref{eqn:main1_discr}) and
its system counterpart, in section~\ref{sec:matrixform} we
reformulate some of these methods in a matrix-oriented setting,
that explicitly exploits the Kronecker form of the linear part of
the problem. Among these, we consider the implicit-explicit Euler
(IMEX Euler) and 2SBDF methods (IMEX-2SBDF), and low order
exponential integrators (Exp Euler). {Algorithmic details
that make the matrix-oriented methods particularly efficient are
discussed in section~\ref{sec:details}, where the \emph{reduced}
methods rEuler, rExp and rSBDF working in the spectral space are
introduced.}
{In section~\ref{sec:heat} we discuss the application of the
matrix-oriented methods to the semilinear Heat equation, representative of
(\ref{eqn:main1}) as test RD-PDE when the exact solution is known.
We present a stability and convergence {experimental}
study for the proposed
numerical schemes together with comparisons of computational costs
in order to emphasize the advantages of solving sequences of
{matrix problems with respect} to the usual vector approach.}
{In section~\ref{sec:sys}, we extend the above reduced
schemes to deal with RD-ODE matrix systems \eqref{eqn:matrixform_sys}.
In section~\ref{sec:models}, we present the features of
\emph{diffusion-driven instability} or \emph{Turing theory}
\cite{been,Murray} for pattern formation in nonlinear RD-PDE
systems and the challenges for the numerical approximation of Turing
pattern solutions; we solve the prototype Schnackenberg model \cite{Murray} and the DIB-morphochemical model
\cite{EJAM2015}, representative of these pattern formations in different application contexts.
All reported experiments are performed in Matlab \cite{matlab} on a quadcore processor Intel Core(TM) i7-4770 CPI@3.40GHz, 16Gb RAM.

\section{Classical vector methods}\label{sec:classical}
For the time stepping of \eqref{eqn:main1_discr} we can consider
the following methods, where for the sake of simplicity we
consider a constant timestep $h_t >0$ {and the time grid $t_n
= n h_t$, $n=0, 1, \dots, N_t$ so that $(\u_n)_{ij} \approx
u(x_i,y_j,t_n)$ in each point $(x_i,y_j)$ of the discretized space}:
\begin{enumerate}
\item {\it IMEX methods}.

\begin{enumerate}
\item[i)] {First order Euler: We discretized in time as
$\u_{n+1}-\u_n=h_t(A\u_{n+1}+f(\u_{n}))$, so that
\begin{equation}\label{eqn:IMEXvector}
 (I-h_t A) \u_{n+1}= \u_n+h_t f(\u_n), \quad n=0, \dots, N_t-1,
\end{equation}}
where $\u_0 $ is given by the initial condition in
\eqref{eqn:main1_discr}; the linear part is treated implicitly, while
the reaction (nonlinear) part $f$ is treated explicitly
\cite{Ascher95,Frank97,Ruuth95}.

\item[ii)] Second order SBDF.
The widely used IMEX 2-SBDF method \cite{Ruuth95,Ascher95} applied
to \eqref{eqn:main1_discr} yields
\begin{equation}
3 \bu_{n+2} - 4\bu_{n+1} + \bu_{n} = 2 h_t A \bu_{n+2} + 2 h_t (2 f(\bu_{n+1}) - f(\bu_{n})), \quad n= 0, 1, \dots , N_t -2\\
\label{eqn:2-SBDF}\end{equation} As usual, $\bu_0$ is known, while
a step of the first order IMEX-Euler scheme can be used to determine
$\bu_1$ (\cite{Ruuth95,Ascher95}).

\end{enumerate}

\item {\it Exponential integrator.} Exponential {first order} Euler 
method \cite{hoc}:
\begin{eqnarray}\label{eqn:exp}
\u_{n+1} = e^{h_t A} \u_n + h_t \varphi_1(h_tA) f(\u_n)
\end{eqnarray}
where $e^{h_t A}$ is the matrix exponential, 
and $\varphi_1(z) = (e^z -1)/z$ is the first ``phi'' function \cite{hoc}.
\item {\it ADI method}.  We consider the two-stage time
stepping {when ${\ell}(u)=\Delta u = u_{xx}+u_{yy}$ is the
Laplace operator}
that approximates the reaction term in explicit way:
\begin{equation}\label{ADIij}
\begin{split}
& \frac{u_{ij}^{n+\frac 1 2}-u_{ij}^n}{h_t/2}=(u_{xx})_{ij}^{n+\frac 1 2}+(u_{yy})_{ij}^{n}+f(u_{ij}^{n}), \\
& \frac{u_{ij}^{n+1}-u_{ij}^{n+\frac 1
2}}{h_t/2}=(u_{xx})_{ij}^{n+\frac 1
2}+(u_{yy})_{ij}^{n+1}+f(u_{ij}^{n}).
\end{split}\end{equation}
Let $U_n\approx U(t_n) \in\RR^{N_x\times N_y}$. After discretization we obtain 
\begin{equation}\label{ADI}
\begin{split}
&\left(I-\frac{h_t}{2} T_1\right)U_{n+\frac 1 2}=
\left( I+ \frac{h_t}{2} T_1\right)U_n + \frac{h_t}{ 2} F(U_n) \\
&{U_{n+1}\left(I- \frac{h_t}{2} T_2^T\right)= U_{n+\frac
1 2} \left( I+ \frac{h_t}{2} T_2^T\right) + \frac{ h_t}{ 2}
F(U_n)}.
\end{split}
\end{equation}
We remark that the ADI method naturally treats the approximation in matrix terms, therefore it is the 
closest to our methodology.

\end{enumerate}

\section{Matrix formulation of classical methods}\label{sec:matrixform}
In this section we reformulate some of the time steppings of
section~\ref{sec:classical} in matrix terms, by exploiting the
Kronecker sum in (\ref{eqn:kron}). We then provide implementation
details to make the new algorithms more efficient. {To this
end, we recall that $U(t)$ defines the matrix whose elements
approximate at the time $t$ the values of $u(x_i,y_j,t)$ at the
nodes $(x_i, y_j)$ in the discrete space}, {and $U_n\approx U(t_n)$ such that $\u_n=vec(U_n)$}. We shall see that the
matrix-oriented approach leads to the evaluation of matrix
functions and to the solution of linear matrix equations with
small matrices, instead of the solution of very large vector
linear systems. We stress that the {matrix formulation}
does not affect the convergence and stability properties of the
underlying time discretization method. Rather, it exploits the
structure of the linear part of the operator to make the
computation more affordable. {In particular, this allows one to
refine the space discretization, so as to capture possible
peculiarities of the problem; see, e.g., section~\ref{sec:Schnakenberg}.}

In the following we derive the time iteration associated with the
single differential equation (\ref{eqn:main1}) yielding the semi-discrete ODE matrix system  (\ref{eqn:matrixform}) . A completely analogous iteration will be obtained for the ODE matrix system
(\ref{eqn:matrixform_sys}), see section~\ref{sec:sys}.

\subsection{Matrix-oriented first and second order IMEX methods}\label{sec:IMEX_method}

The matrix-oriented versions of the IMEX methods rely on the Kronecker form of
$A$ in (\ref{eqn:kron}) and on its property that allows to trasform the vector linear system
into a matrix linear equation to be solved, of much smaller size.

{Consider the discretized times $t_n = n h_t, \ n=0, \dots N_t$ with
timestep $ h_t >0$.}
{Then adapting the one-step scheme in (\ref{eqn:IMEXvector}),  to the differential
matrix form \eqref{eqn:matrixform}, yields

$$
U_{n+1} - U_n = h_t (T_1 U_{n+1} +U_{n+1}T_2) + h_t F(U_n),
$$
which, after reordering, gives the following linear matrix equation, called the {\em Sylvester} equation,
\begin{eqnarray}\label{eqn:imex_mat}
(I-h_t T_1)U_{n+1} +U_{n+1}(-h_t T_2) = U_n + h_t F(U_n),  \quad
n=0, \dots, N_t-1  .
\end{eqnarray}
Therefore, to obtain the next iterate $U_{n+1}$ the matrix approach for the IMEX Euler method requires the
solution of a Sylvester equation at each time step, with coefficient
matrices $(I-h_t T_1)$, $(-h_t T_2)$ and right-hand side $U_n + h_t F(U_n)$.
The numerical solution of \eqref{eqn:imex_mat} is described in section~\ref{sec:details}. 
The matrix equation \eqref{eqn:imex_mat} should be compared with the vector form, requiring the
solution of a linear system of size $N_xN_y\times N_xN_y$ at each
time step. It is important to realize that for a two-dimensional
problem on a rectangular grid, the number of nodes required in
each direction need not exceed a thousand, even in the case a fine
grid is desired to capture possibly pathological behaviors. Hence,
while the Sylvester equation above deals with, say, matrices of
size $500\times 500$, the vector form deals with matrices and
working vectors of size $250\,000\times 250\,000$. Arguably, these
latter large matrices are very sparse and structured, so that
strategies for sparse matrices can be exploited; nonetheless, the
Sylvester equation framework allows one to employ explicit
factorizations, also exploiting the fact that the coefficient matrices do not
change with the time steps. Algorithmic details will be given in
section~\ref{sec:details}.

As second order IMEX strategy for (\ref{eqn:matrixform}) we consider the two-step
method {IMEX-2SBDF} seen in \eqref{eqn:2-SBDF} for the vector formulation.
For the matrix form, given the initial condition $U_0$, and a further approximation $U_1$
-- obtained for instance by the IMEX Euler method -- at each time step $t_{n+2}$ the method
determines the following matrix equation
$$
3U_{n+2}-4U_{n+1} + U_n = 2h_t\left( T_1U_{n+2} + U_{n+2}T_2 + 2 F(U_{n+1})-F(U_n)\right ) ,
$$
which, after reordering, leads
once again to the solution of a Sylvester equation in the unknown matrix $U_{n+2}$ given by
$$
 \left (3I-2h_t T_1\right) U_{n+2}+U_{n+2}\left (-2h_t T_2\right ) =4U_{n+1}- U_n +2h_t (2F(U_{n+1})-F(U_n)),  \quad
n=0, \dots, N_t-2 .
$$
The coefficient matrices are
$3I-2h_tT_1$, $(-2h_tT_2)$ and the right-hand side is $4U_{n+1}- U_n +2h_t (2F(U_{n+1})-F(U_n))$.

\subsection{Exponential Euler method}\label{sec:exp}
A matrix-oriented version of the exponential Euler approach can
exploit (\ref{eqn:kron}) in the computation of both the
exponential and the phi-function. In particular, the following
property of the exponential matrix is crucial
$$
e^{h_t A}=e^{h_t(I \otimes T_1+T_2^T\otimes I)}= e^{h_t T_2^T}\otimes
e^{h_tT_1} . 
$$
Therefore, for $\u={\rm vec}(U)$ we have
$$
e^{h_t A} \u = \left(e^{h_t T_2^T}\otimes e^{h_tT_1}\right)  \u  
= {\rm vec}(e^{h_tT_1} U e^{h_t T_2}) .
$$
Moreover, the operation $ v = h_t \varphi_1(h_t A) f =
A^{-1}(e^{h_t A}f - f)$ can be performed by means of a two step
procedure which, given $F$ such that $f={\rm vec}(F)$ delivers $V$
such that $v={\rm vec}(V)$:
\begin{enumerate}
\item[-] Compute  $\ G =   e^{h_tT_1} F e^{h_t T_2}$ 
\item[-] Solve $\ T_1 V + V T_2 = G - F \quad \quad $ for $V$.
\end{enumerate}
 Therefore, the matrix-oriented version of the Exponential Euler method first computes
the matrix exponential of multiples of $T_1$ and $T_2$ once for all. Then, it
obtains the approximation $U_{n+1}$ by solving a Sylvester matrix equation at each time step. More precisely,
\begin{enumerate}
\item Compute  $E_1 = e^{h_t T_1}$,  $E_2 = e^{h_t T_2^T}$ 
\item For each $n$
\begin{eqnarray}
{\rm Solve} &\quad T_1 V_n + V_n T_2 = E_1F(U_n) E_2^T - F(U_n)  \label{eqn:sylvexp}\\
{\rm Compute} & \quad U_{n+1} = E_1 U_n E_2^T + V_n .\nonumber
\end{eqnarray}
\end{enumerate}
Several  implementation suggestions are given in the
next section. 
It is important to realize that to be able to solve (\ref{eqn:sylvexp}) the two
matrices $T_1$ and $-T_2$ must have disjoint spectra. Unfortunately, Neumann boundary
conditions imply that both $T_1$ and $T_2$ are singular, leading to a zero
common eigenvalue. To cope with this problem we employed the following
differential matrix equation, mathematically equivalent  to (\ref{eqn:matrixform}),
\begin{eqnarray}\label{eqn:flex_meq}
\dot U = (T_1 - \sigma I) U + U T_2 + (F(U) + \sigma U).
\end{eqnarray}
{With this simple ``relaxation'' procedure
the matrix $T_1-\sigma I$ is no longer singular, and has
no common eigenvalues with $-T_2$, at the small price of including an extra linear term
to the nonlinear part of the equation.
We note that adding and subtracting the term $\sigma U$ to the ODE may be beneficial -- though
not strictly necessary --
also for the other methods; thus in section~\ref{sec:stability}
we include a stability analysis for all considered time integration strategies
based on the relaxed matrix equation (\ref{eqn:flex_meq}). }
\subsection{Implementation details}\label{sec:details}
Whenever the matrix sizes are not too large, say up to a thousand,
the previously described matrix methods can be made more efficient by
computing a-priori a spectral decomposition of the coefficient matrices involving
$T_1$ and $T_2$. In the following we shall assume that the two matrices
are diagonalizable, so that their eigenvalue decompositions can be determined.
Let them be $T_k = X_k \Lambda_k X_k^{-1}$, $k=1,2$, with $X_k$ nonsingular and
$\Lambda_k={\rm diag}(\lambda_1^{(k)}, \lambda_2^{(k)}, \ldots)$ diagonal.

Let us first consider the IMEX Euler iteration in (\ref{eqn:imex_mat}).
Compute the
$N_x\times N_y$ matrix $L_{i,j} = 1/((1-h_t\lambda_i^{(1)})+(-h_t\lambda_j^{(2)}))$. Hence, at each iteration $n$
we can proceed as follows
\begin{itemize}
\item[1.] Compute $\widehat U_n = X_1^{-1} Q(U_n) X_2 \quad$ where $ \quad Q(U_n) = U_n + h_t F(U_n)$;
\item[2.] Compute $U_{n+1} = X_1 (L\circ \widehat U_n) X_2^{-1}$
\end{itemize}
where $\circ$ is the Hadamard (element by element) product. The second step performs
the solution of the Sylvester equation by determining the solution
entries one at the time, in the eigenvector bases, {and then the result is
projected back onto the original space to get $U_{n+1}$} \cite{Simoncini.survey.16}.
Proceeding in the same manner, the corresponding version for
IMEX-2SBDF can be derived. Letting this time $L_{i,j} =
1/((3-2h_t\lambda_i^{(1)})+(-2h_t\lambda_j^{(2)}))$ at each
time iteration $n$ we have:
\begin{itemize}
\item[1.] Compute $\widehat U_n = X_1^{-1}Q(U_{n}, U_{n+1})X_2 \quad$ where $ \quad Q(U_{n}, U_{n+1}) = 4U_{n+1}- U_n +2h_t
(2F(U_{n+1})-F(U_n))$;
\item[2.] Compute $U_{n+1} = X_1 (L\circ \widehat U_n) X_2^{-1}$.
\end{itemize}
In the following numerical experiments we
will call these methods: {\bf reduced IMEX-Euler (rEuler)} and
{\bf reduced 2SBDF (rSBDF)}. Whenever the RD-PDE problem is linear,
that is $f(u) = \alpha u + \beta$, the computation further
simplifies, since all time steps can be performed in the
eigenvector basis, and only at the final time of
integration the approximate solution is interpolated back to the
physical basis; see section~\ref{sec:heat}.

In a similar way, the matrix-oriented exponential
Euler integrator described in section~\ref{sec:exp} can be
rewritten as
\begin{enumerate}
\item Compute  $\widehat e_k = {\rm diag}(e^{h_t \lambda_1^{(k)}}, e^{ht\lambda_2^{k)}}, \ldots)$,
$k=1,2;$ $\widehat E=\widehat e_1 \widehat e_2^*$ and
$\widehat L_{i,j} = (h_t\lambda_i^{(1)}+h_t\lambda_j^{(2)})^{-1}$, with $\widehat E, \widehat L \in{\mathbb C}^{N_x\times N_y}$.
\item For each $n$,
\item[] $\hskip 0.2in$ Compute $\widehat F_n = X_1^{-1} F(U_n) X_2$   \hskip 0.4in \% Project $F(U_n)$ on the eigenbases;
\item[] $\hskip 0.2in$ Compute $G =  \widehat E \circ \widehat F_n - \widehat F_n$ \hskip 0.55in \% Apply {\tt exp} and form the Sylvester eqn rhs;
\item[] $\hskip 0.2in$ Compute $V = \widehat L \circ G$ \hskip 0.9in \% Solve the Sylvester eqn;
\item[] $\hskip 0.2in$ Compute $U_{n+1} = X_1 (\widehat E \circ (X_1^{-1} U_n X_2) + V)  X_2^{-1} $ \hskip 0.2in \% Compute the next iterate
\end{enumerate}
In the following numerical
experiments we will call this method {\bf reduced Exp (rExp)}.

\noindent {If the ``relaxation'' approach corresponding to \eqref{eqn:flex_meq} is considered, 
the quantity $F_\sigma(U) =  F(U) + \sigma U$ replaces $F(U)$ in the algorithm above, while
the spectral decomposition of $T_1(\sigma)=T_1-\sigma I$ replaces that of $T_1$.}

\section{The semilinear heat equation}\label{sec:heat}

We start by specializing the matrix methods of the previous sections to the
case of the following Heat Equation (HE)
with linear source (reaction) term and zero Neumann boundary
conditions (BCs):
\begin{equation}\begin{cases} u_t=d \Delta u + \alpha u \quad (x,y)\in \Omega \subset \mathbb{R}^2,\ t \in ]0, T], \ \\
(\boldsymbol n \nabla u)_{|\partial \Omega}=0, \quad
u(x,y,0)=u_0(x,y), \label{Heat1}\end{cases}\end{equation} with
$d\in \mathbb{R_+}$ the diffusion coefficient, $\alpha \in
\mathbb{R}$ the reaction coefficient and $\Omega =[0, \ell_x] \times
[0, \ell_y]$ a rectangular domain.
We consider the initial condition $u_0(x,y)=A_0 \cos(c_x x)\cos(c_y y)$,
for which the exact solution of \eqref{Heat1} is given by $u^*(x,y,t)=
e^{(\alpha-(c_x^2+c_y^2)d)t}u_0(x,y)$.

For the numerical treatment, we consider a finite
difference approximation for spatial derivatives based on the
{Extended Central Difference Formulas} (ECDF) \cite{amo1,TuringJCAM2012}.
These schemes consider the approximation of the Neumann BCs with
the same order of schemes used in the interior domain, so that no
reduction of order arises near the boundaries. In particular, we
apply the scheme of order $p=2$ as follows.
 Let us discretize the domain $\Omega$ with $N_x$ and $N_y$
 interior points, giving step sizes $h_x=\ell_x/(N_x+1)$ and
$h_y=\ell_y/(N_y+1)$. Let $U_{ij}(t)\approx u(x_i,y_j,t)$, for
$i=1,\dots, N_x$, $j=1,\dots,N_y$ be the values of the approximate solution at
the interior mesh nodes, and let $\bu = {\rm vec}(U)$.

Let $T_x \in \mathbb{R}^{N_x\times N_x}$ and $T_y \in
\mathbb{R}^{N_y\times N_y}$ be the usual tridiagonal matrices
corresponding to the approximation of the second order derivatives
in \eqref{Heat1} by central differences (order $p=2$), along the
$x$ and $y$ directions, and  zero Neumann BCs approximation.  More precisely,
$T_x=diag(1,-2,1)+B$, and similarly for $T_y$, with corresponding dimensions, 
where the BCs term (see \cite{TuringJCAM2012,JCAM2016}) is given by
\begin{equation} \label{eqn: Bneu}
B=\frac{2}{3}\begin{bmatrix}
2 & -\frac{1}{2} & \cdots & 0 & 0\\
0 & 0 & \cdots & \cdots & 0\\
\vdots & & & & \vdots \\
0 & \cdots & &-\frac{1}{2} & 2
\end{bmatrix}  .
\end{equation}

Therefore, the semi-discretization of \eqref{Heat1} in vector form is given by
\begin{equation} \label{VectorFormHE}
\dot \u= A \u+ \alpha \u \quad  \u(0)=\u_0,
\end{equation}
where
\begin{equation}\label{eqn:A}
 A= d \widetilde \Delta, \qquad
\widetilde \Delta=\frac{1}{h_x^2}(I_{y} \otimes T_x)+\frac{1}{h_y^2}(T_y \otimes I_{x})
\in \mathbb{R}^{N_xN_y \times N_xN_y}
\end{equation}
while its matrix counterpart becomes
\begin{equation} \label{MatrixForm}
\dot U= T_1 U+ UT_2+ \alpha U, \quad  \quad  U(0)=U_0, \quad
T_1=\frac{d}{h_x^2} T_x, \quad T_2=\frac{d}{h_y^2} T_y^T .
\end{equation}
The solution of the Heat equation can take full advantage of the matrix formulation
{because of the linearity of the reaction term, see section \ref{sec:details}}.
In Algorithm \ref{Alg1} we report the reduced IMEX-Euler method all in
the eigenvector space. Note that the computational cost is thus kept to the minimum.

\begin{algorithm}
\caption{Reduced-Sylvester method for the Heat Equation} \label{Alg1} 
\begin{algorithmic}
\STATE Compute $X_k, \Lambda_k$, $k=1,2$ such that $T_1=X_1 \Lambda_1 X_1^{-1}$, $T_2=X_2 \Lambda_2 X_2^{-1}$
\STATE $\hat U_1=X_1^{-1}U_0 X_2$
\STATE $L_{i,j}=1./(\Lambda_1(i,i) + \bar\Lambda_2(j,j))$
\FOR{$n=1:N_t$}
\STATE $\hat U_1=L \circ( \hat U_1+\alpha h_t \hat U_1)$
\ENDFOR
\STATE $U_1=X_1 \hat U_1 X_2^{-1}$
\end{algorithmic}
\end{algorithm}

\noindent If time discretization is performed by the reduced Euler exponential integrator {\bf rExp} 
we can take full advantage of the linearity of the operator, by
including the constant term $\alpha$ into the $x$-direction
matrix, thus avoiding the use of the term involving the function
$\varphi_1$. {The resulting expression for the
next iterate is thus given by
\begin{equation} \label{HE_omega1}
\bu_{n+1} = e^{h_t
(A + \alpha I)} \bu_n, \end{equation}
and $h_t (A + \alpha I) =
I\otimes (T_1+h_t \alpha I) + T_2^T\otimes I$. By first computing
the eigenvalue decompositions of $T_1+ h_t \alpha I$ and of $T_2$
(which give the corresponding matrices $E$, $X_1$ and $X_2$), we
obtain the matrix iteration in the eigenvector basis,
$$
\widehat U_0= X_1^{-1} U_0 X_2, \qquad \widehat U_{n+1} = \widehat E \circ \widehat U_n , \,\, n=0,1,\ldots
$$
where the elements of $\widehat E$ are all possible products of
eigenvalues of $T_1+ h_t a I$ and $T_2$ (see the general description in section~\ref{sec:details}).}
Only at the final
integration time the solution is projected back onto the physical basis
as $U_* = X_1 \widehat U_* X_2^{-1}$.

\subsection{Stability analysis}\label{sec:stability} In this
subsection, we present a stability analysis of the schemes
proposed in the previous sections, by considering as PDE test
problem the Heat Equation in \eqref{Heat1} and its vector form
semi-discretization \eqref{VectorFormHE}. By using the spectral
decomposition $ A = Q \Lambda Q^{-1}$, we can define $\widetilde
\u = Q^{-1} \u $, and work with the associated {scalar test
problem} for each component $\widetilde u =(\widetilde \u)_j,
\lambda= \lambda_j$ where $\lambda_j = (\Lambda)_{jj}$. Note that
it holds that $\lambda = d \widetilde \lambda$ where $ \widetilde
\lambda_j$ is an eigenvalue of $\widetilde \Delta$ in
(\ref{eqn:A}).\footnote{Since both matrices 
$T_x$ and $T_y$ are similar to symmetric matrices, their
eigenvalues are real; moreover, it can be shown that for the
considered matrix $B$, their eigenvalues are non-positive.}
In general, we can thus
consider the scalar ODE test problem
\begin{equation}\label{ScalarTest}
 \dot{\widetilde u} = \lambda \widetilde u + \alpha \, \widetilde u
\end{equation}
where $\lambda \in \RR_-$ models the diffusion and
$\alpha \in \RR$ models the reaction.
The exact solution of
\eqref{ScalarTest} is $u^*(t) = u_0 e^{(\lambda + \alpha)t}$ that
goes to zero for $t \rightarrow \infty$ (i.e. it is asymptotically stable) if and only if
\begin{equation}\label{eq1}
 \lambda + \alpha < 0 \quad \Longleftrightarrow \quad  \alpha < -\lambda.
\end{equation}
If $h_t$ is the time
stepsize, let be $\xi = \lambda h_t$ and $\mu = \alpha h_t$, then the region in the plane $(\xi,\mu)$, with $\xi
< 0$ where stationary solutions can be obtained is given by the
half-plane $\mu < -\xi$.
(In the following with some abuse of notation we use $u$ instead of
$\widetilde u$ in \eqref{ScalarTest}.)
{Following (\ref{eqn:flex_meq}), we relax the equation with}
$\sigma u = -\omega \alpha u$ in \eqref{ScalarTest} as follows:
\begin{equation}\label{ScalarTestw}
\dot u = (\lambda u -\sigma u) + (\alpha u+\sigma u) \  \Longleftrightarrow \ \dot u = \lambda_\omega u + \alpha_\omega \,u, \quad \quad   0 \leq
\omega \leq 1
\end{equation}
where
$\lambda_\omega =(\lambda +\omega \alpha)$ and
$\alpha_\omega=(1-\omega) \alpha$.
For $\omega \in [0,1]$ the problem \eqref{ScalarTestw} is equivalent to
\eqref{ScalarTest}, that is it admits the same solution $u^*$, and for
$\omega=0$ we find exactly \eqref{ScalarTest}. Moreover, by
\eqref{eq1} we need to enforce the constraints
\begin{eqnarray}\label{eqn:constr}
\lambda_\omega < 0 , \qquad
\lambda_\omega +\alpha_\omega = \lambda +\alpha < 0.
\end{eqnarray}
Our aim is twofold: (i) identify the stability
regions $R_{met}$ in the $(\xi,\mu)$ plane of the considered methods,
where these schemes are able
to reproduce the asymptotic stationary behavior of the
theoretical solution; (ii) identify the restrictions on the timestep $h_t$ 
in terms of the diffusion and reaction terms.

We recall that the relaxed approach needs to be used to apply the
rExp integrator to the general nonlinear model, since
$A$ in (\ref{VectorFormHE}) is singular. Although in this linear case the general Euler exponential
method would not be required, we believe that the analysis still provides
valuable indications of the stability properties of the methods towards the
general nonlinear setting.

\noindent The {\bf IMEX Euler method} for \eqref{ScalarTestw} 
 is given by $u_{n+1}= u_n + h_t ( \lambda_\omega u_{n+1} + \alpha_\omega u_n)$, from which we get
\begin{equation}\label{EIMEXw}
(1-\lambda_\omega  h_t) u_{n+1} = (1+ \alpha_\omega  h_t) u_n.
\end{equation}
Letting $\xi = \lambda_\omega h_t < 0 $ and $\mu =\alpha_\omega
h_t$, we have $ \displaystyle u_{n+1} = \frac{(1+ \mu)}{(1-\xi)}
u_n . $ Then, taking into account the constraints in
(\ref{eqn:constr}), the EIMEX numerical solution is
asymptotically stable in (see Figure~\ref{FigStabReg} (left) dashed lines)
$$
{\cal R}_{EIMEX} = \left\{ (\xi,\mu) \in \RR_- \times \RR \ | \quad
 \displaystyle \frac{|1+\mu|}{|1-\xi|} < 1 \right\} \quad
\Leftrightarrow \quad    \xi - 2 < \mu < -\xi .
$$

\begin{figure}[!t]
\centering
\includegraphics[height=6cm, width=7.95cm]{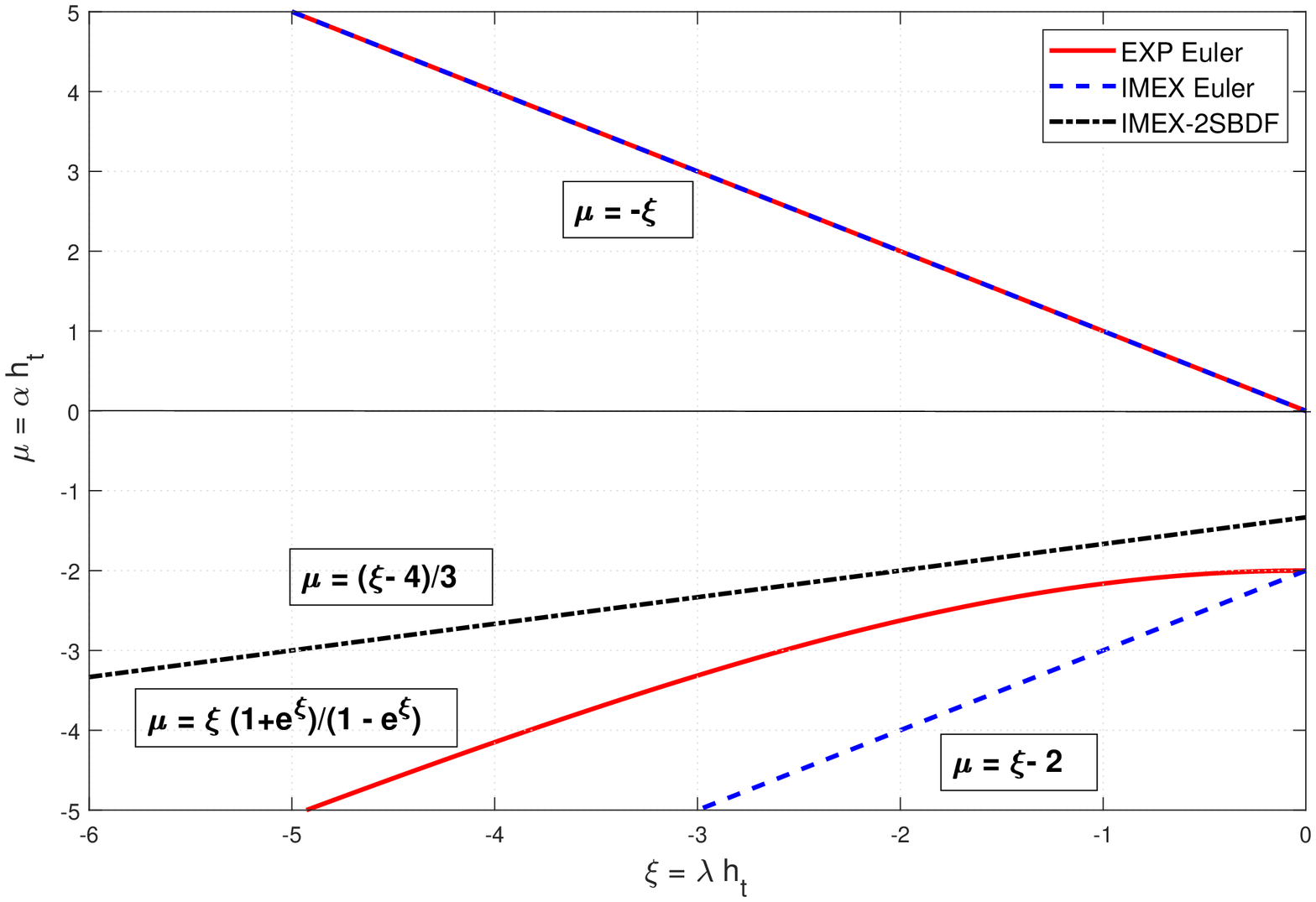} 
\includegraphics[height=6cm, width=7.95cm]{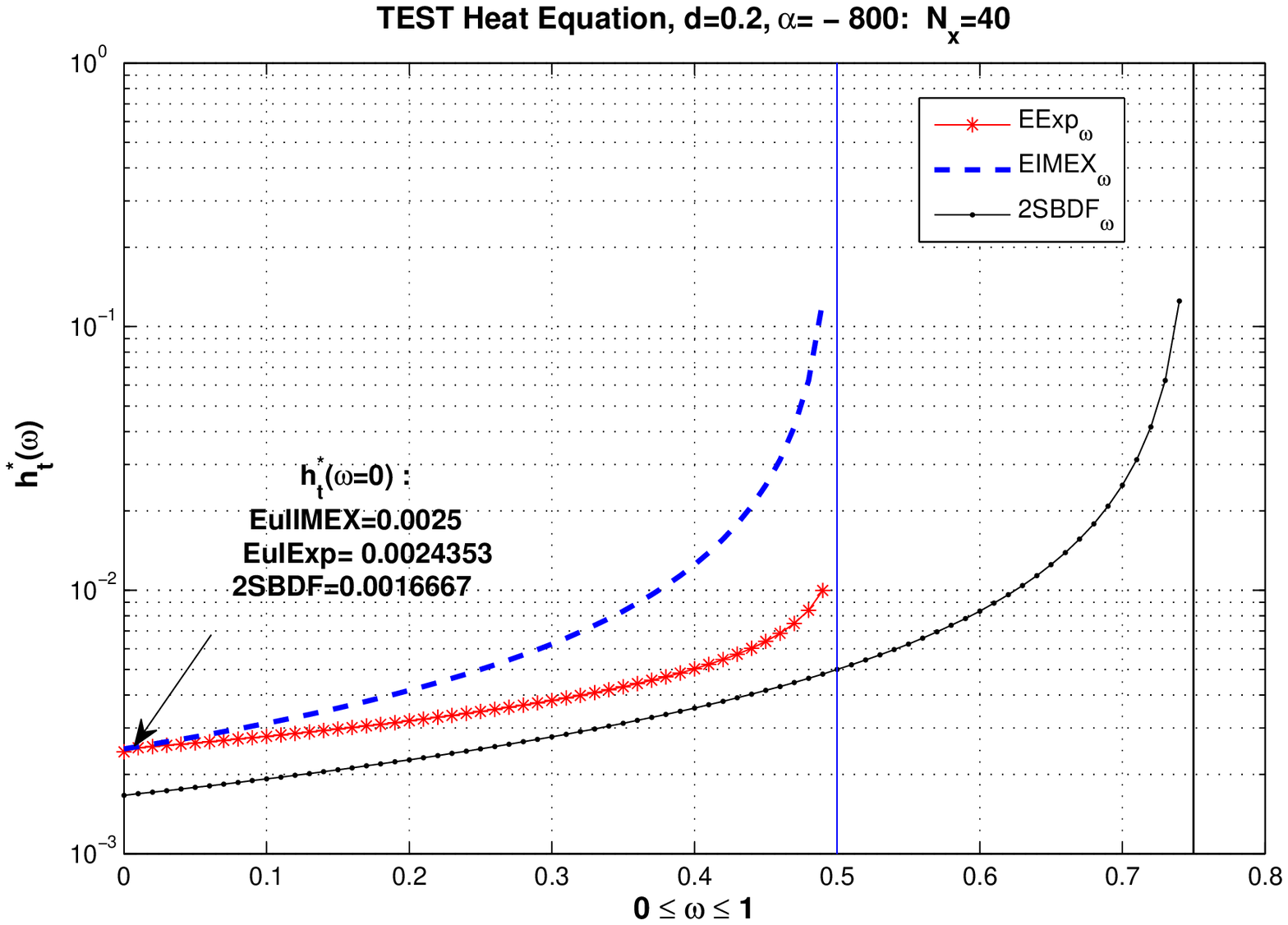}
\caption{Heat equation - Left plot: stability regions in the plane $(\xi,\mu)$
for IMEX and EXP Euler methods and the 2-SBDF method applied to
the linear test problem \eqref{VectorFormHE}.  Right plot: For $d=0.2$, $\alpha=-800$ in \eqref{Heat1} and $N_x=N_y=40$ meshpoints in the spatial domain $\Omega =[0,1] \times [0,1]$, we show the critical timesteps of all methods in dependence of the relaxation parameter $\omega \in [0, 1]$.}
\label{FigStabReg}
\end{figure}

Explicitly writing $\mu$ and $\xi$, the bound $ \mu > \xi -2$ allows us to
determine possible timestep restrictions. Indeed, the
bound $ \mu > \xi -2$ corresponds to
\begin{equation}\label{eq2}
{\cal F}_{EIMEX}(h_t,\lambda,\alpha,\omega):=[(1-2\omega)\alpha -\lambda]h_t + 2 > 0 .
\end{equation}

A detailed analysis shows that this condition is always satisfied for all considered
parameters, except if $(1-2\omega)\alpha < \lambda < 0$, in which case
the timestep constraint
$ h_t < 2/(\lambda-(1-2\omega)\alpha)$. 
Since this condition must hold for all $\lambda$s, we can determine the
critical stepsize in correspondence of the worst case as 
\begin{equation}\label{htIMEX}
 h_t < \frac{2}{\lambda_M -(1-2\omega)\alpha} =: h_{crit}(\omega) ,
\end{equation}
where  $\lambda_M= \max_{j=1,\dots, N_x \cdot N_y} \lambda_j(A) = d \max_j \widetilde
\lambda_j(\widetilde \Delta) < 0$. The curve of the critical stepsize $h_{crit}(\omega)$, as a function of
$\omega \in [0,1]$, is shown in Figure~\ref{FigStabReg} (right) for the
case $d=0.2$, $\alpha=-800$, and
$\lambda_M$ computed numerically for $N_x=N_y=40$ meshpoints
in the spatial domain $\Omega =[0,1] \times [0,1]$. 
{It is easy to see that the timestep restriction arises only for $\omega < 1/2$}.

\noindent A similar analysis for the {\bf Exponential Euler method}
determines the following stability region (see Figure~\ref{FigStabReg} (left) red continuous line):
$$
{\cal R}_{EEXP} = \left\{ (\xi,\mu) \in \RR_- \times \RR \ | \quad
\displaystyle| \mu \frac{e^\xi -1}{\xi} + e^\xi \displaystyle | < 1 \right\} \quad
\Leftrightarrow \quad    \xi \frac{ 1 + e^\xi}{1 - e^\xi} < \mu < -\xi .
$$
By substituting $\xi = \lambda_\omega h_t$ and $\mu =
\alpha_\omega h_t$ the lower bound for $\mu$ implies the timestep
restriction $-(\lambda_\omega + \alpha_\omega) e^{\lambda_\omega h_t} - 
(\lambda_\omega - \alpha_\omega) >0$, which is equivalent to requiring
\begin{equation}\label{eq3}
 {\cal F}_{EEXP}(h_t,\lambda,\alpha,\omega):=
-(\lambda + \alpha) e^{(\lambda +\omega \alpha) h_t} +[(1-2\omega)\alpha-\lambda] >0 .
\end{equation}
Once again, a detailed analysis of the sign shows that only for
$(1-2\omega)\alpha < \lambda < 0$ we determine a constraint  by
identifying the critical timestep $h_t(\omega)$ as $\omega$ varies.
This is obtained by looking for the zeros $z^*$ of
${\cal F}_{EEXP}(z,\lambda_M,\alpha,\omega)$ for $\omega \in [0,1]$. For the
same PDE data as for the IMEX Euler method, the curve of these zeros
$z^*=h_t^*(\omega), \omega \in [0,1]$ is reported in
Figure~\ref{FigStabReg} (right).  
{Again, the timestep restriction arises only for $\omega < 1/2$}.

\noindent As shown in \cite{TuringJCAM2012}, the stability region for 
 the second order {\bf IMEX 2-SBDF method} (\ref{eqn:2-SBDF})
in the
half-plane $(\xi, \mu)$, $\xi < 0$ can be studied by the roots
$|z_{1,2}(\xi,\mu) | $ of the second order characteristic
polynomial associated to (\ref{eqn:2-SBDF}) when applied to the
linear test problem:$
(3 - 2 \xi ) z^2 - 4(\mu +1)z + 1+ 2 \mu =0.
$
Hence, we have the stability region given by (see Figure~\ref{FigStabReg} (left)  dash-dot line)
$$
 {\cal R}_{SBDF} = \left\{(\xi,\mu) \in \RR_- \times \RR \ | \quad
|z_{1,2}(\xi,\mu) | \leq 1 \right\} \quad \Leftrightarrow \quad
\frac{ \xi -4}{3} < \mu < -\xi.
$$
Simple algebra shows that the same
result is obtained if we consider the relaxed test problem
\eqref{ScalarTestw} and $\xi = \lambda_\omega h_t$ and $\mu =
\alpha_\omega h_t$.
Again, the lower bound for $\mu$ implies
the restriction $(3\alpha_\omega -\lambda_\omega) h_t < -4$ on the timestep $h_t$, 
which is equivalent to requiring
\begin{equation}\label{eq4}
 {\cal F}_{SBDF}(h_t,\lambda,\alpha,\omega)=[(3-4\omega)\alpha -\lambda]h_t + 4 > 0
\end{equation}
Similar arguments as above imply that the constraint on $h_t$ is given
if $ \ (3-4\omega)\alpha < \lambda < 0$, in which case $ h_t <
\displaystyle \frac{4}{\lambda-(3-4\omega)\alpha}$.
Having to hold for all $\lambda$s, the
critical timestep can be found by the worst case as:
\begin{equation}\label{htSBDF}
 h_t < \displaystyle \frac{4}{\lambda_M -(3-4\omega)\alpha} =: h_{crit}(\omega) .
\end{equation}

In Figure~\ref{FigStabReg} (right) we report the critical timesteps
\eqref{htSBDF} $h_t^*(\omega)$, for $\omega \in [0,1]$, for again the same PDE data.
In this case $ \alpha < 0 $ and stepsize restrictions arise only for
$ 0 \leq \omega < 3/4$. 

For the case $\omega=0$, that is without relaxation,
the critical timesteps for the considered schemes are:
$h_t^{EIMEX} = 0.0025, h_t^{EEXP} =0.0024353,h_t^{2SBDF} =
0.0016667$. (Note that these values can change for different space discretization values of $N_x,N_y$.)
By increasing the relaxation parameter $\omega$ until $\omega^*$
less stringent bounds arise, that is for all schemes an
improvement in the stability requirements are obtained. The best
gain is obtained by EIMEX, followed by EEXP and then by the more
demanding second order scheme. This behavior
is to be expected, due to the size of the respective stability regions
shown in Figure~\ref{FigStabReg} (left). For $\omega \geq \omega^*$ the
methods are unconditionally stable, with $\omega^* = 0.5$ for
EEIMEX and EEXP, with $\omega^*=0.75$ for 2SBDF. It is worth
noting that the application of the reduced Exponential Euler
method for the Heat Equation in \eqref{HE_omega1} corresponds to
the relaxation $\omega=1$ and then, for $\alpha < 0$, the method
will not require any timestep restriction, then it is
unconditionally stable.

\begin{remark}
For the scalar linear ODE \eqref{ScalarTest} 
the relaxed $\omega-$schemes correspond to applying in time 
an implicit approximation
for the diffusion part in $\lambda$, and a classical $\theta-$ method for the reaction part.
Hence, the relaxed IMEX Euler \eqref{EIMEXw} for $\omega=0$ is equivalent
to $EIMEX$, whereas for $\omega=1$ it corresponds to the fully implicit Euler method.
The determined stability properties are thus expected.  For the other relaxed schemes, we are not aware of similar stability results.
It is worth mentioning that whenever the nonlinear term $f(u)$ includes a linear term with negative coefficient, 
the $\sigma$-correction in \eqref{eqn:flex_meq} simply corresponds to moving this linear term 
to the diffusion part of the expression.
\end{remark}

\subsection{Numerical results}
In this section we experimentally explore the performance of the
considered methods, that is
{\bf rEuler}, {\bf rsbdf}, {\bf ADI} and {\bf rExp}, on the simple model
matrix problem 
\begin{equation} \label{eqn:Heatmatrix}
\dot U = T_1 U+ U T_2 + \alpha U, \quad U(0) =U_0.
\end{equation}
where $T_1,T_2$ are given in \eqref{MatrixForm}.
According to (\ref{eqn:flex_meq}), for rExp we actually solve the differential matrix equation $\dot U = (T_1+\alpha I) U+ U T_2$.

\noindent Throughout this set of experiments we solve \eqref{Heat1} on $\Omega =[0,1] \times [0,1]$, with the parameter choice $\alpha=-800$,
$d=0.2$,  $c_x=2$, $c_y=1$, $A_0=1$, $T_f=0.1$ and $N_x=N_y$ in the spatial meshgrid. 
\begin{figure}[t]
\centering
\includegraphics[width=8cm, height=6cm]{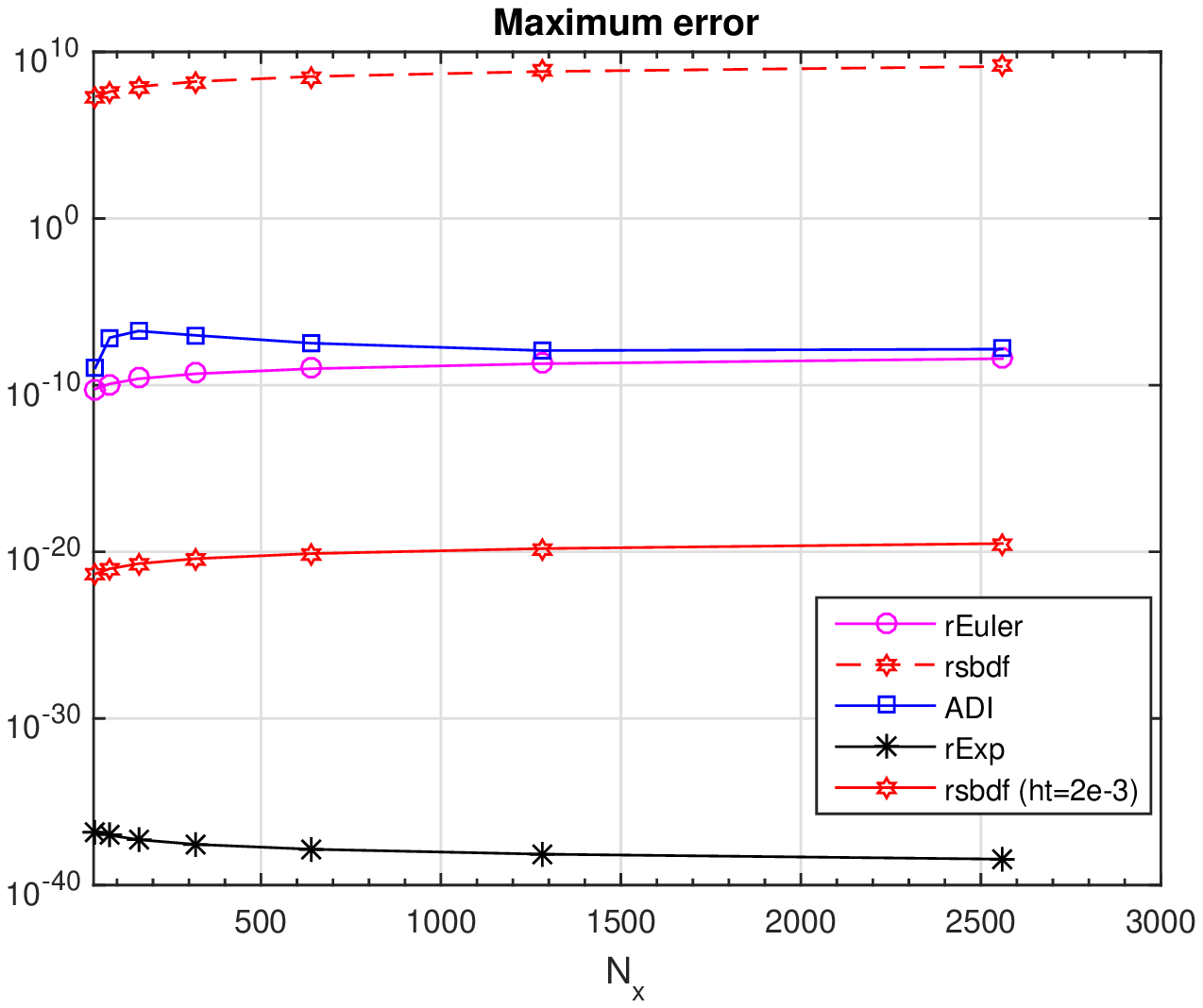}%
\includegraphics[width=8cm, height=6cm]{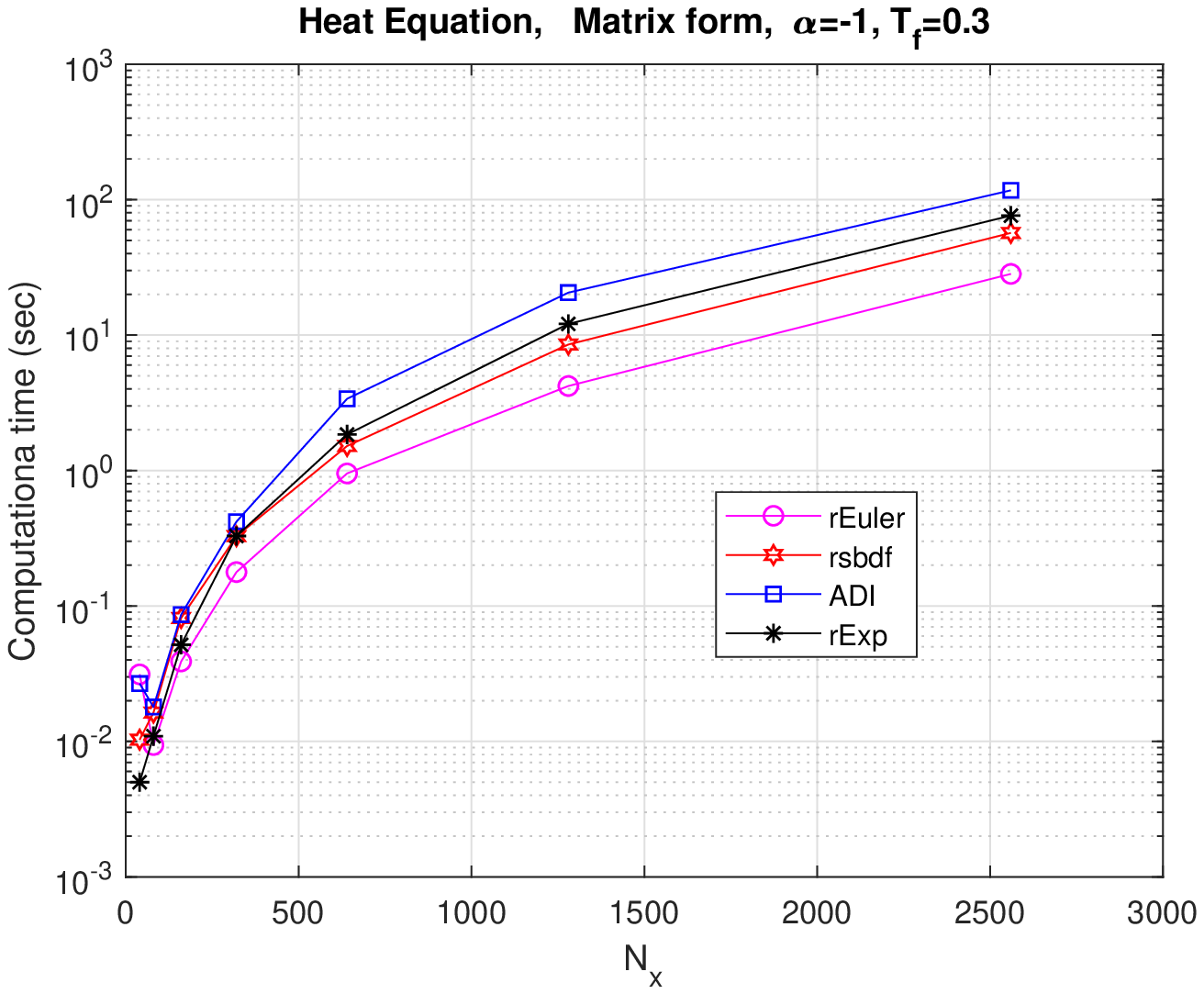} 
\caption{Heat equation with $d=0.2$, $\alpha=-800$. Left plot: 
Error of various \emph{reduced} time integrators in matrix form for $h_t=10^{-3}$ 
(and also for $h_t=2\, 10^{-3}$ in rSBDF). Right plot: 
Computational costs for all methods for $h_t=10^{-3}$.
\label{fig:stability}}
\end{figure}

To experimentally test the stability analysis of the previous section in the
matrix case as the problem dimension increases, we consider two critical time steps, and vary $N_x$ as 
$N_x^k=40 \cdot 2^k$ for $k=0,... ,6$. 
{The left plot of Figure~\ref{fig:stability} reports
the behavior of the maximum (spatial) error for the solution at the final 
time $T_f$ for all methods for $h_t=10^{-3}$ (for rSBDF the
history for $h_t=2\, 10^{-3}$ is also reported).
The plot confirms the analysis of the previous section:
for the value $h_t=2\, 10^{-3}$ falling outside the region of
stability of the method, rSBDF provides an unacceptably high error, whereas 
the error behaves as expected by the convergence theory
for the smaller value $h_t=1\, 10^{-3} < h_t^{cr}$.} 
{The plot also confirms} that for the linear problem the exponential method
is exact in exact arithmetic, thus the displayed error is only due to 
the use of finite precision arithmetic in the computation of the spectral information.

Figure \ref{fig:stability}(right) reports on the computational cost (time in seconds) of all
methods as the problem size grows. For the considered matrix dimensions, 
all methods behave {somewhat similarly, although the use of
the reduced strategy appears to be very beneficial}. 
Note that for the largest values of $N_x$, the
Kronecker form of the problem would have a diffusion matrix of size $6\cdot 10^6$,
which would be extremely hard to handle for a large number of times steps, while
limiting the type of time stepping methods to be employed.
These considerations are particularly important in the case of systems of nonlinear
equations, as considered in the next section.
To illustrate this point, in Table~\ref{tab:Heat} we report the
CPU times of a full run with the standard vector-oriented IMEX scheme, compared
with the reduced Euler approach. 
For this experiment, we consider $\alpha=-1$ and $T_f=0.3$, so that there is 
no restriction on the time step $h_t$. Moreover, we fix the ratio 
$\nu=\frac{h_t}{h_x^2}$ to preserve the order of convergence. 
Starting from $h_{t}^1=10^{-2}$, $N_x=32$ and $h_{x}=\frac{1}{N_x+1}$, 
$h_t$ is halved at each step $k$, that is $h_t^k:=\frac{h_{t}}{2^k}$ 
for $k=0,1,2...$, and $h_x^k:=\frac{h_x^{k-1}}{\sqrt{2}}$.
In the vector case, an LU factorization with pivoting
of the coefficient matrix is performed a-priori\footnote{To increase sparsity
of the factors, {\tt symamd} reordering
was performed on the given matrix.}, so that only sparse triangular
solves are performed at each iteration with the $N_x^2\times N_x^2$ factors.
For larger values of $N_x$, a preconditioned iterative solver
would be required to solve such very large linear system.
\begin{table}[htb]
\begin{center}
\begin{tabular}[]{|c|r|r|r|r|r|r|r|r|r|}
\hline
$N_x$ & $32$ & $45$  & 64 &     90  &128        &180     &256  & 362   \\
\hline
\hline
 Vector &        0.0414  &  0.0115 &    0.0398 &    0.1366 &    0.7899 &    3.1548 &   13.750 &   66.562 \\
\hline
rEuler &   2.3e-3   & 2.7e-3  & 6.3e-3  &  0.0143   & 0.0336 &   0.0878&    0.1824&    0.5699 \\
\hline
\end{tabular}
\end{center}
\caption{Heat equation for $\alpha=-1$ and $d=0.2, T_f=0.3$: Computational 
times in seconds for the 
IMEX Euler method solved in vector form by direct solver and in matrix form by the reduced method 
rEuler.\label{tab:Heat}}
\end{table}

\section{Reaction-diffusion PDE systems: matrix approach}\label{sec:sys}
In this section we apply the  matrix-oriented approach 
to an RD-PDE model with nonlinear reaction terms and zero 
Neumann boundary conditions, given in general by:
\begin{equation}\label{reacdiff}
\begin{cases} u_t =d_1 \Delta u + f_1(u,v), 
\quad (x,y)\in \Omega \subset \mathbb{R}^2,\ t \in ]0, T] \\
v_t= d_2 \Delta v+ f_2(u,v),\\
(n\nabla u)_{|\partial \Omega}=(n\nabla v)_{|\partial \Omega}=0\\
u(x,y,0)=u_0(x,y), v(x,y,0)=v_0(x,y)
\end{cases}
\end{equation}

The Method of Lines for \eqref{reacdiff},
 based on the finite difference space discretization in \eqref{MatrixForm}, yields the 
following system of ODE matrix equations:
\begin{equation}\label{reacdiff_matrix}
\begin{cases} U' = d_1(T_1U + UT_2) + F_1(U,V)\\
V'= d_2(T_1 \, V + V \, T_2) + F_2(U,V)\\
U(0)=U_0, V(0)=V_0 .
\end{cases}
\end{equation}
As in section \ref{sec:matrixform}, the matrix form of the classical ODE methods can be derived for
\eqref{reacdiff_matrix}, giving rise to 
 the solution of the following Sylvester matrix equations 
at each timestep $t_n$,
\begin{equation}\label{sys:one_step}
\begin{cases}
S_1U_{n+1}+U_{n+1}S_2=Q_1^n,\\
R_1 V_{n+1}+V_{n+1}R_2=Q_2^n , \quad n=0, \dots, N_t-1 \quad
\quad U_0, V_0 \ \mbox{given}
\end{cases}
\end{equation}
where $Q_k^n=Q_k^n(U_n, V_n), k=1,2$ in the case of a one step method 
({like IMEX Euler method in the previous sections}) and
$Q_k^n=Q_k^n(U_{n-1}, V_{n-1},U_n,V_n)$ (with $U_0, U_1$ given) for a two-step scheme.
Recalling the procedure of section~\ref{sec:IMEX_method},
for IMEX-Euler we have
\begin{eqnarray}
S_1&=&I-h_t d_1 T_1, \,\, S_2=-h_t d_1 T_2, \quad  R_1=I-h_td_2 T_1, \,\, R_2=-h_t d_2 T_2, 
 \label{rEuler}
\\
Q_1^n&=&U_n+h_tF_1(U_n,V_n), \quad \, Q_2^n=V_n+h_tF_2(U_n,V_n),\nonumber
\end{eqnarray}
while for IMEX-2SBDF we have
\begin{eqnarray} 
S_1&=&3I-2h_t d_1 T_1,\,\, S_2=-2h_t d_1 T_2,\,  \quad  
R_1=3I-2h_td_2 T_1,\, R_2=-2h_t d_2 T_2, \nonumber \\
Q_1^n&=&4U_{n}-U_{n-1}+2h_t(F_1(U_{n},V_{n})-F_1(U_{n-1},V_{n-1})),  
\label{rSBDF}
\\
Q_2^n&=&4V_{n}-V_{n-1}+2h_t(F_2(U_{n},V_{n})-F_2(U_{n-1},V_{n-1})).\nonumber
\end{eqnarray}
Analogously to section~\ref{sec:exp}, also the matrix-oriented version of the
exponential method can be derived for the RD systems. In particular, letting once
again $E_{1,1} = e^{h_t d_1 T_1}$,  $E_{1,2} = e^{h_t d_1 T_2^T}$,
and $E_{2,1} = e^{h_t d_2 T_1}$,  $E_{2,2} = e^{h_t d_2 T_2^T}$
we obtain
{
\begin{eqnarray}\label{rExp}
U_{n+1} &=& E_{11} U_n E_{12}^T + Y_n, \,\, \quad {\rm where} \quad\, (d_1 T_1-\sigma I) Y_n + Y_n (d_1 T_2) = E_{11} \widetilde F_1(U_n,V_n) E_{12}^T \\
V_{n+1} &=& E_{21} V_n E_{22}^T + Z_n, \,\, \quad {\rm where} \quad \, (d_2T_1-\sigma I) Z_n + Z_n (d_2 T_2) = E_{21}\widetilde F_2(U_n,V_n) E_{22}^T ,
\end{eqnarray}
where $\sigma$ is as described in section \ref{sec:exp},
while $\widetilde F_1(U_n,V_n) = F_1(U_n,V_n) +\sigma U_n$ 
and $\widetilde F_2(U_n,V_n) = F_2(U_n,V_n) +\sigma V_n$.}
In particular, the approach requires the solution of two Sylvester 
equations per step, which is the same cost as
for the IMEX procedure, together with matrix-matrix multiplications 
with the exponentials. As already discussed
for the single equation case, these costs can be significantly 
reduced by working in the eigenvector basis
of $T_1$ and $T_2$. 
A Matlab implementation of the rEuler and rExp methods is reported in the Appendix.

\section{Nonlinear reaction-diffusion systems with Turing solutions}\label{sec:models}

Reaction-diffusion systems like \eqref{reacdiff} arise in several
scientific applications and describe mathematical models where the
unknown variables can have different physical meaning, for example
chemical concentrations, cell densities, predator-prey population
sizes, etc. Ranging from ecology to bio-medicine, depending on
the parameters in the model kinetics $f_1, f_2$, different kinds of
solutions can be studied, for example traveling waves or
oscillating dynamics. We are interested in
\emph{diffusion-driven or Turing instability} solutions of \eqref{reacdiff},
arising from the perturbation of
a stable spatially homogeneous solution. More precisely, let $(u_e,v_e)$ be the stable
solution to the homogeneous equations \eqref{reacdiff} with no
diffusion, that is $f_1(u_e,v_e)=0=f_2(u_e,v_e)$. {Adding diffusion can force
spatial
instability to take place}, leading to an asymptotic interesting spatial
pattern (so-called {\it Turing} pattern),
characterized by structures like spots, worms, labyrinths, etc. More recently, in
\cite{Neubert97,Neubert2002} the authors proved that the transient
dynamics is important for pattern formation. In particular, the
concept of \emph{reactivity} describing the short-term transient
behavior, is necessary for Turing instabilities.
Let $\w =(u,v)$ and let
$
J = J(u_e,v_e) = \begin{bmatrix}
f_{1,u} & f_{1,v} \\
f_{2,u} & f_{2,v}
\end{bmatrix}|_{\w_e}
$
be the Jacobian of the linearized ODE system associated to
\eqref{reacdiff} evaluated at the spatially homogeneous solution
$\w_e =(u_e,v_e)$.
The spatially homogeneous solution $\w_e =(u_e,v_e)$ is stable 
if the eigenvalues of $J$ all have negative real part, that
is $J$ is a stable matrix.
In \cite{Neubert97} the authors defined $\w_e$ as
\emph{reactive equilibrium} if the largest eigenvalue of the
symmetric part of $J$ is positive:
$$ 
\lambda_{\max}(H(J)) > 0, \quad \quad H(J) = (J+J^T)/2.
$$
We recall that $J$ may be stable, while $H(J)$ is {\it not}
stable, that is $H(J)$ has at least one strictly positive
eigenvalue. This key feature is typical of stable highly non-normal matrices
\cite{Trefethen.Embree.05}.
If the initial condition $\w_0(x,y)=(u_0(x,y),v_0(x,y))$ in
\eqref{reacdiff} is a small (random) perturbation to $\w_e$, the
RD-PDE solution in the initial transient, say $\v(x,y,t)$, is
governed by the linearization
$$
\v_t = D \, \Delta \v + J \v, \quad D=diag(d_1, d_2).
$$
By applying the Fourier transform $ \widetilde \v (\k, t) =
\int_{-\infty}^{\infty} e^{i (k_x x+ k_y y) t} \v(x,y,t) dx dy$ \
this equation becomes the linear ODE system
$$
\widetilde \v' = \widetilde J \widetilde \v, \quad \quad
\widetilde J = J - \| \k \|_2^2 D,$$ 
where $\k =(k_x,k_y)$
and $\| \k \|_2^2= k_x^2+k_y^2= (\pi \nu_x /\ell_x)^2+ (\pi \nu_y
/\ell_y)^2$ accounts for the spatial frequencies $\nu_x,\nu_y$. 
The Turing theory ensures that if the largest real part of the
eigenvalues of $\widetilde J$ is positive for some $\k$, then
perturbations with this spatial frequency will grow in time and produce
spatial patterns, so that $\w_e$ is destabilized by diffusion. The
Turing conditions on the model parameters identify a range of
spatial modes such that the pattern formation arises for 
$\|\k\|_2^2 \in [\k_1^2, \k_2^2]$ (see e.g. \cite{been}). In
\cite{Neubert2002}, the authors show that the largest eigenvalue
of $H(J)$ must be positive for an eigenvalue of
$\widetilde J$ to have positive real part. {\it Reactivity} is
therefore a prerequisite for pattern formation via Turing
instability. It would be desiderable that numerical methods for
the approximation of Turing patterns also reproduced
the reactivity features during the initial transient regime, in
addition to reaching  the asymptotic stability.
%
 The numerical approximation of Turing pattern solutions is thus
challenging for the three following reasons: (i) longtime integration is
needed to identify the final pattern as asymptotic solution of the PDE
system; (ii) the time solver should capture
the reactivity phase at short times;
(iii) a large finely discretized space domain $\Omega$ is
required to carefully identify the spatial structures of the
Turing pattern.
The matrix-based procedures described in the previous sections
allow us to efficiently address all three items above in
the numerical treatment of typical models for this physical
phenomenon.
In particular, our numerical tests will highlight items (i) and (ii)
with the well-known Schnackenberg model
\cite{been,Murray}. We will apply the ADI method \eqref{ADI} often
used in the literature (see, e.g., \cite{TuringJCAM2012}) and the Sylvester
based methods studied in the previous sections
\emph{rEuler},\emph{rExp} and \emph{rSBDF}. 
To deal with item (iii), we
propose to apply the matrix-based approach to solve the
morpho-chemical model (briefly said DIB model) recently proposed
in the literature to study the morphology and the chemical
distribution in an electrodeposition process typical of
charge-richarge processes in batteries. 

\begin{figure}[hbt]
\centering
\includegraphics[width=2.5in, height=2.5in]{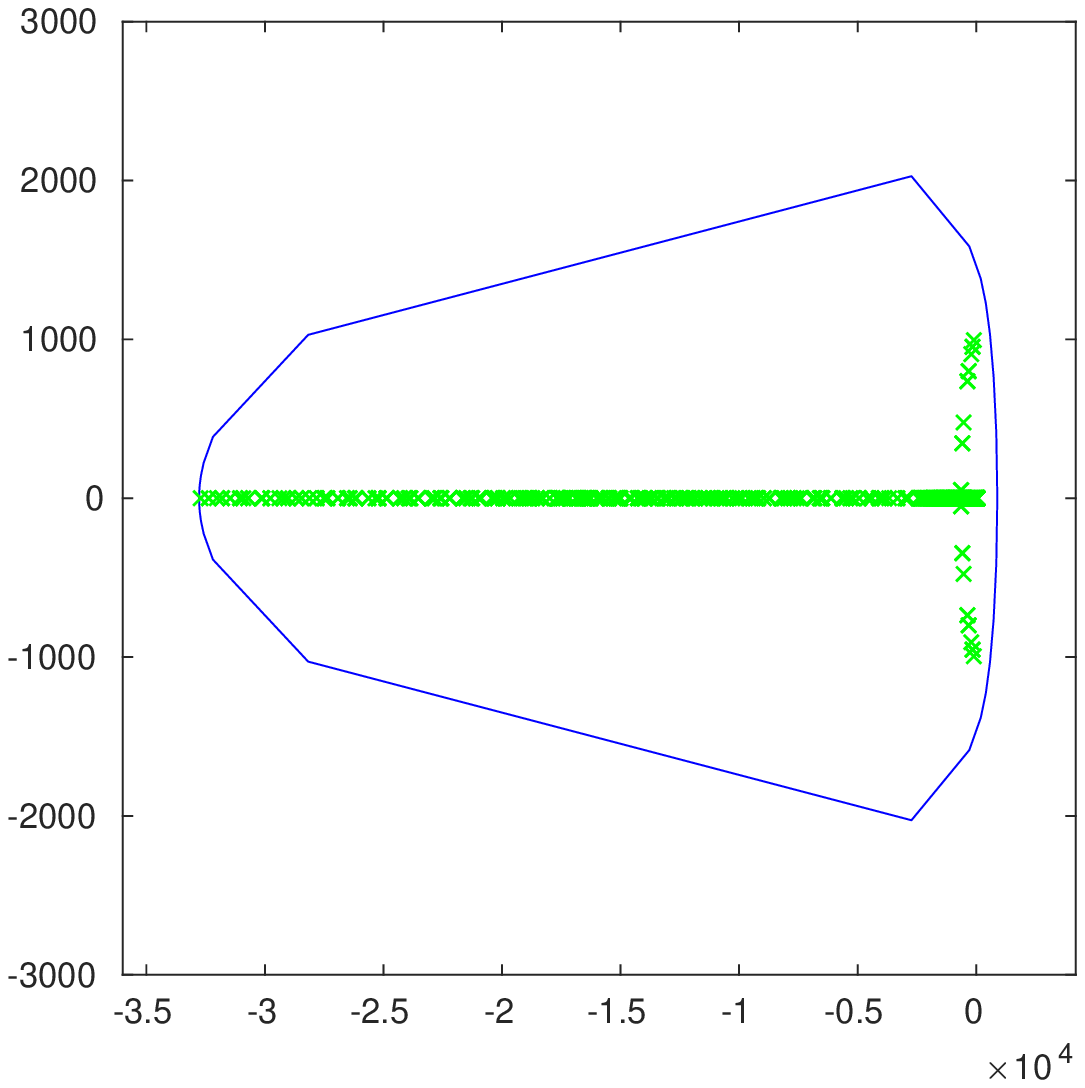} \hskip 0.8in
\includegraphics[width=2.5in, height=2.5in]{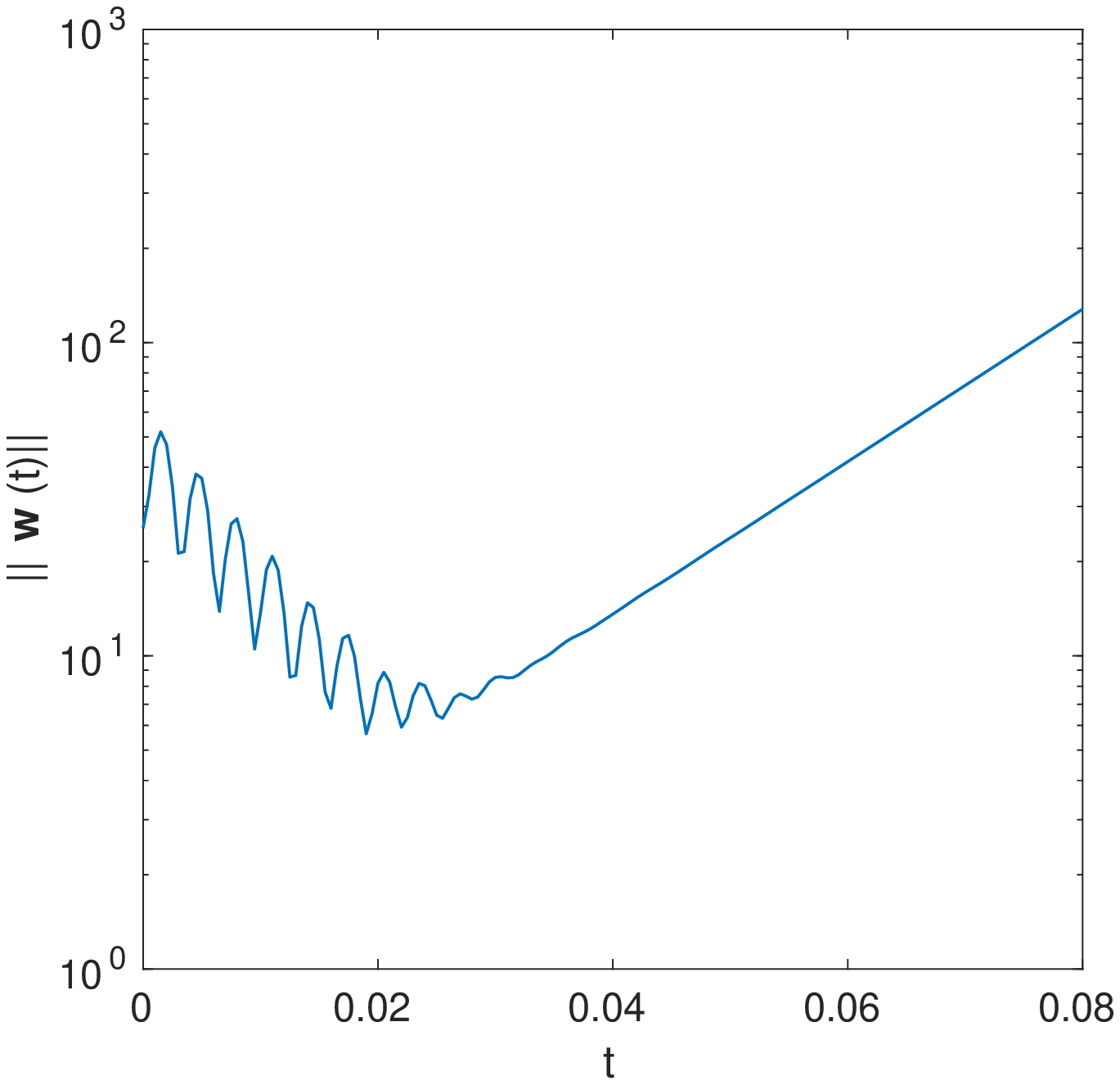}
\caption{Schnackenberg model. Left plot: Eigenvalues and Field of values of $\widetilde J$.
Right plot: $\|w\|$ with $w=e^{t \widetilde J} w_0$.  \label{fig:fov}}
\end{figure}

\subsection{Schnakenberg model}\label{sec:Schnakenberg}
The RD-PDE for the Schnakenberg model is given by
\begin{equation}\label{Schnak1}
\begin{cases} u_t =\Delta u + \gamma(a - u+ u^2 v), \quad (x,y)\in \Omega= [0, 1] \times [0, 1],\ t \in ]0, T_f] \\
v_t= d\Delta v+ \gamma(b - u^2 v),\\
(n\nabla u)_{|\partial \Omega}=(n\nabla v)_{|\partial \Omega}=0\\
u(x,y,0)=u_0(x,y), v(x,y,0)=v_0(x,y)
\end{cases}\end{equation}
This model has received great attention in the recent
literature (see, e.g., \cite{liu,ric}) because in spite of its
simplicity, it is representative of classical patterns typically found
in biological experiments. The model parameters
$a$, $b$, $d$, $\gamma$ are positive constants, and a unique
stable equilibrium exists, which undergoes the Turing instability,
given by $u_e=a+b$, $v_e= \frac{b}{(a+b)^2}$. We consider the
typical values $d=10, \gamma=1000, a=0.1, b=0.9$,
yielding a cos-like spotty pattern of the type $\cos(\nu_x \pi
x)\cos(\nu_y \pi y)$ with the selected modes
$(\nu_x,\nu_y)=(3,5),(5,3)$ \cite{been} (see 
Figure~\ref{fig:SchnackSolCost}). We consider the initial conditions
$u_0(x,y)=u_e+10^{-5}{\tt rand}(x,y), v_0(x,y)=v_e+10^{-5}{\tt rand}(x,y)$
where {\tt rand} is the default Matlab function with fixed seed of the
generator (rng('default')) at the beginning of each simulation.
To numerically confirm the spectral analysis of the previous
section, we consider the matrices (see (\ref{eqn:A}))
$$
\widetilde J = \begin{bmatrix} \tilde\Delta & 0 \\ 0 & A \end{bmatrix} +  J , \quad \quad J=\gamma
\begin{bmatrix}
(-1 + 2 u_e v_e) {\mathbb I} & u_e^2 {\mathbb I} \\
-2 u_e {\mathbb I}  &  - u_e^2 {\mathbb I}
\end{bmatrix}  .
$$

Here ${\mathbb I}$ is the matrix of all ones. Note that
$J$ has (multiple) eigenvalues $\lambda(J)_{\pm}= \gamma(-0.1 \pm 0.99499i)$, while
$\lambda(H(J))_{1,2}\in\{ -1.0849\gamma, 0.88489\gamma\}$.
The equilibrium is
reactive because the largest eigenvalue of $H(J)$ is positive. For
the chosen parameters, the left plot of Figure~\ref{fig:fov}
displays the field of values and eigenvalues of $\widetilde J$, from which we
can see that the matrix $\widetilde J$ is not stable (its largest (real) eigenvalue
is about 60), and its symmetric part is not negative
definite\footnote{The eigenvalues of $H(\widetilde J)$ are contained in
the interval $\RR\cap W(\widetilde J)$, where the field of values 
$W(\widetilde J)$ of an $N\times N$ matrix $\widetilde J$ is defined as
$W(\widetilde J) = \{ z \in \CC, z=(x^*\widetilde Jx)/(x^*x), x\in\CC^N\}$.}
The right plot of Figure~\ref{fig:fov} indicates
norm divergence of the solution to the linearized problem, $w = e^{t \widetilde J}
w_0$, thus confirming the existence of the reactivity phase in the 
initial transient regime. {This behavior has been observed
for instance in} population dynamics ODE models in \cite{Trefethen.Embree.05}, 
{however we are not aware of a similar  computational evidence 
of the reactivity concept  for Turing pattern formation. 
The experimental results in Figure~\ref{fig:fov} thus support  
the reactivity analysis in \cite{Neubert97,Neubert2002} for the
numerical treatment of the Schnakenberg model.}

\begin{figure}[htb]
\centering
\includegraphics[scale=0.52]{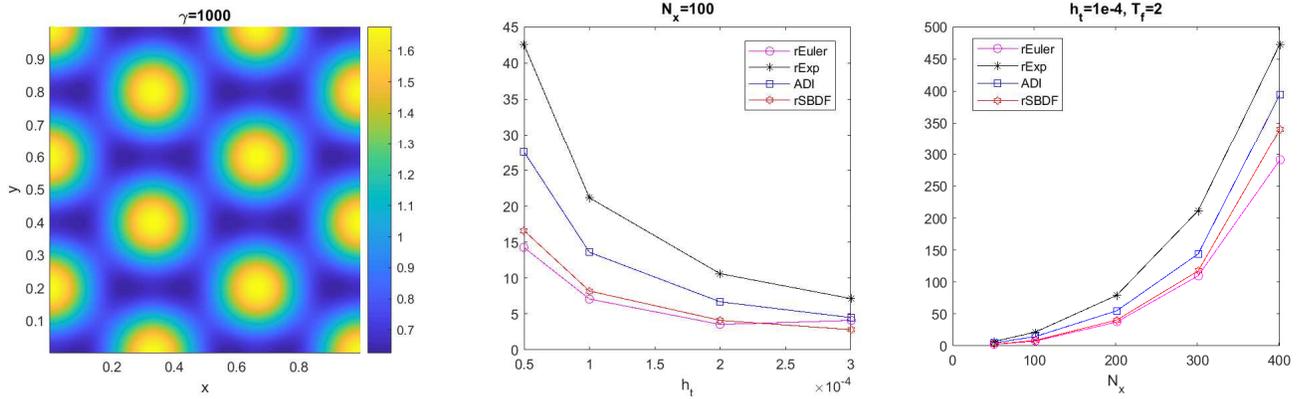}
\caption{Schnackenberg model.  Left plot: Turing pattern 
solution for $\gamma=1000$ ($N_x=400$). Center plot: CPU times 
(sec) for Test (a), $N_x=100$ variation of $h_t$. 
Right plot: CPU times (sec) for Test (b), $h_t =10^{-4}$, increasing 
values of $N_x=50,100,200,300,400$.}
\label{fig:SchnackSolCost}
\end{figure}

\noindent To study the time dynamics in our simulations we will report the
values of the space mean value
\begin{equation} \label{mean}
\langle U_n \rangle = {\rm mean}(U_n) 
\approx \langle u(t_n) \rangle = \frac{1}{|\Omega|} 
\int_\Omega u(x,y,t_n)\, dx\, dy  \quad \quad t_n=n \, h_t, \quad n=0,\dots , N_t ,
\end{equation}
that for $t \rightarrow T_f = N_t \, h_t$ will tend to a constant
value, say $ \langle U_n \rangle \rightarrow {\bar u}$, if a
stationary pattern is attained. We also report the behavior of
the increment $ \delta_n =\| U_{n+1}-U_n \|_F$ (Frobenius norm)
that will tend to zero if the steady state is reached. These two
indicators will be useful also to describe the numerical behaviors
of the methods in the initial transient and then to study their
reactivity features. Plots with $V_n$ show a completely analogous
behavior. For the final time $T_f= 2$ and $N_y=N_x$, we
present the following two tests.
\begin{figure}[hbt]
\centering
\includegraphics[scale=0.54]{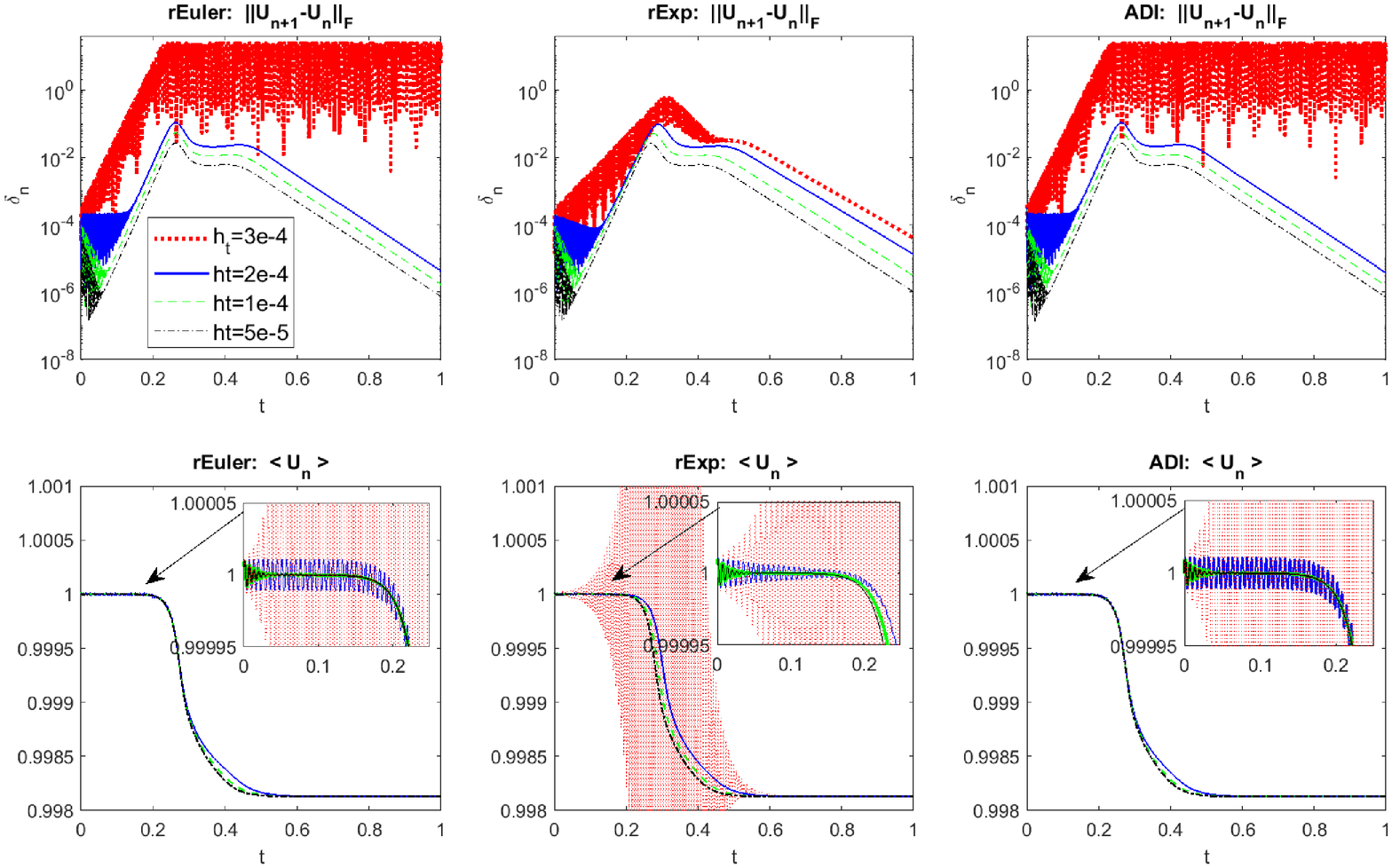}
\caption{Schnackenberg model- Test (a). Indicators for the rEuler, rExp and ADI methods with $N_x=100$ fixed and varying $h_t$ in $\{0.5\cdot 10^4, 10^{-4}, 2\cdot 10^{-4}, 3\cdot 10^{-4}\}$. We show the time behaviors of the increments $\delta_n = \| U_{n+1}-U_n \|_F$ (upper subfigures) and of the space mean values $<U_n>$ (lower subfigures). the zoom insets highlight the reactivity zones ${\cal I}_1$ for each method.}
\label{fig:Test_a}
\end{figure}

\begin{figure}[hbt]
\centering
\includegraphics[scale=0.52]{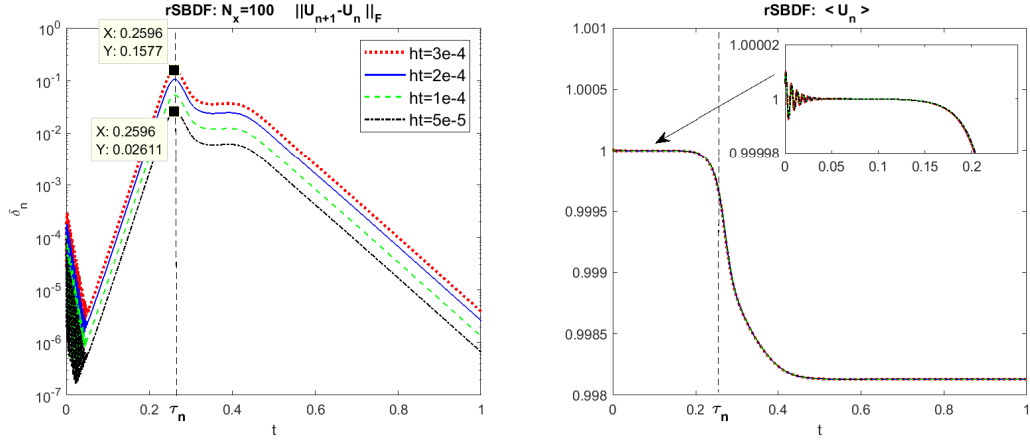}
\caption{Schnackenberg model- Test (a). Indicators for rSBDF $N_x=100$ fixed and varying $h_t$ in $\{0.5\cdot 10^4, 10^{-4}, 2\cdot 10^{-4}, 3\cdot 10^{-4}\}$. Left plot: time behaviors of the increments $\delta_n = \| U_{n+1}-U_n \|_F$, the reactivity zone ${\cal I}_1$ and the stabilizing zone ${\cal I}_2$ are separated by the critical time value $\tau_n$ . Right plot: space mean values $<U_n>$ with zoom inset in the reactivity zone.}
\label{fig:Test_aSBDF}
\end{figure}

{\bf Test (a)}. We fix $N_x=100$ and vary $h_t$ in $\{0.5\cdot 10^4, 
10^{-4}, 2\cdot 10^{-4}, 3\cdot 10^{-4}\}$. The 
simulations reported in Figure~\ref{fig:Test_a} for the rEuler, rExp, ADI methods and in Figure~\ref{fig:Test_aSBDF} for the rSBDF method, show that all methods have a
similar qualitative behavior and that two time regimes
${\cal I}_1 = [0, \tau]$ and ${\cal I}_2 =] \tau, T_f]$ can be
distinguished. In ${\cal I}_1$ reactivity
holds: the oscillating solution departs from the spatially
homogeneous pattern due to the superimposed (small random)
perturbations and becomes unstable, in ${\cal I}_2$ the solution starts
to stabilize towards the steady Turing pattern. Numerically the
value of $\tau$ can be approximated a-posteriori by $\tau_n$, which is the time value
where the maximum of the
increment $\delta_n$ is achieved. Let us discuss in more details the
characteristics of the different methods. To this end, 
the three upper subplots of Figure \ref{fig:Test_a} and the left plot of Figure \ref{fig:Test_aSBDF} for rSBDF show the increment $\delta_n$ of subsequent approximate solutions as time marches. The three lower subplots of Figure \ref{fig:Test_a} and the right plot of Figure \ref{fig:Test_aSBDF} for rSBDF show
the behavior of the mean value $<U_n>$ as time steps proceed.

{{\it ${\cal I}_1$-reactivity zone}. The upper subplots of 
Figure~\ref{fig:Test_a} show that there exists
an initial phase of oscillations whose length is method dependent,
and for $h_t \rightarrow 0$ this length tends to a certain small 
value $\tau_0$. Then for {$\tau_0 \leq t_n \leq \tau$} the solution {\it must} be
unstable, as the necessary condition for Turing instability
requires. Comparing with the mean values $<U_n>$ curve in the lower subplots of 
Figure~\ref{fig:Test_a}, the transfer from ${\cal I}_1$ to
${\cal I}_2$ corresponds to the steep part of the 
curve $\langle U(t) \rangle$ that connects the very short-term 
and the final states of the system; the value of $\tau$ can be related to the
inflection point of this curve. In this experiment it is
possible to note that $\tau_0=\tau_0(h_t^p)$ and 
$\tau=\tau(h_t^p)$. In fact, as also the zoom insets show, 
for $h_t \rightarrow 0$ rEuler and ADI have the same behavior, rExp has
curves $<U_n>$ with different slopes depending on $h_t$. 
Figure~\ref{fig:Test_aSBDF} for the rSBDF method emphasizes that 
this scheme is able to identify the best approximation of $\tau_0$ and $ \tau$ also
for larger value of $h_t$, as it could be expected because it is a 2nd order
method. The zoom inset shows also that the reactivity oscillations are 
kept to the minimum amplitude compared to the other schemes.
}

{\it ${\cal I}_2$-stabilizing zone}. For $h_t < h_t^{cr}$ all
methods reach the asymptotic pattern. Here for $h_t^{cr}\simeq 3\cdot 10^{-4}$
 rEuler and ADI do not attain any pattern ({for clarity,} 
the erratic oscillations
of $<U_n>$ are only reported in the zooms of Figure~\ref{fig:Test_a}), 
while rExp attains the final pattern after a
fully oscillating transient behavior (see the (red) oscillations
in the central subplots of Figure \ref{fig:Test_a}). 
{This value of $h_t$ may be interpreted as} a critical value for the ``reactive stability'' 
of rExp. In the central plot of Figure~\ref{fig:SchnackSolCost} we
report the computational costs of all methods. We recall that by
varying $h_t$ an increasing number $N_t$ of
Sylvester matrix equations in \eqref{reacdiff_matrix} of the same
dimension need be solved.
{rEuler and rSBDF have almost the same cost and
are cheaper than the other methods. rExp is more expensive than
ADI. It is worth noting that the IMEX schemes in vector form are 
{are known to commonly be more expensive than
the ADI method for this problem} (see e.g. \cite{JCAM2016}).}

\begin{figure}[thb]
\centering
\includegraphics[scale=0.52]{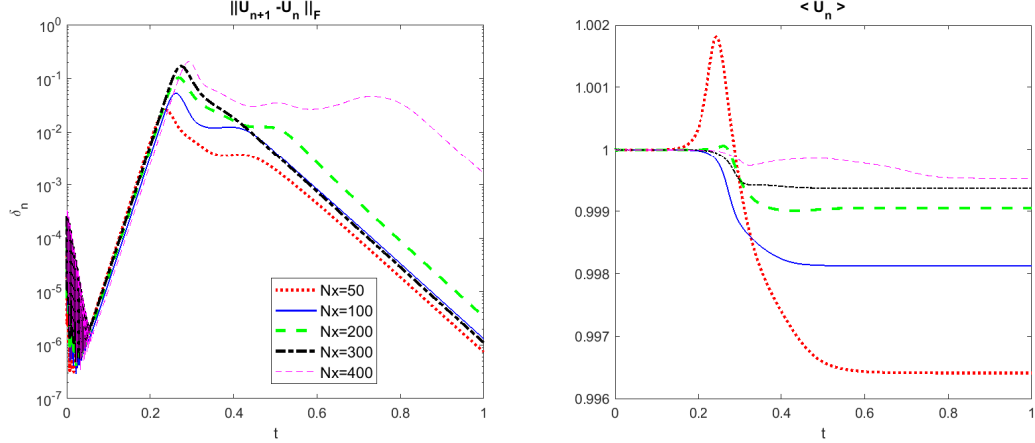}
\caption{Schnackenberg model-Test (b), for rSBDF with $h_t =1$e-4 and increasing values of
 $N_x=50,100,200,300,400$. Left: Increment $\delta_n = \| U_{n+1}-U_n \|_F$.
Right: Space mean value $<U_n>$.}
\label{fig:Test_b}
\end{figure}
{\bf Test (b)}: We fix $h_t=10^{-4}$ and vary $N_x=50,100,200,300,400$,
such that the matrix methods solve the same number $N_t$ of
Sylvester equations of increasing sizes. As a sample,
results for rSBDF are reported in Figure~\ref{fig:Test_b}. All other methods have a
similar time dynamics behavior. 
{The right plot shows that
the} final value of $<U_n>$, that is $\bar u=\bar u(N_x)$, changes with
$N_x$, as expected. {The left plot seems to indicate that
$\tau=\tau(h_t,N_x)$ is} an increasing function of $N_x$.
This might be related to the fact the (discrete)
initial conditions $U_0, V_0$ have different
sizes and include different (though small) random perturbations.
This sensitivity with respect to the
initial conditions is well known in the
pattern formation literature (see, e.g., \cite{Maini}) and it goes under the name
of \emph{robustness problem}. 
This is the main reason why we do not
propose to apply low rank approximation methods for
\eqref{reacdiff_matrix} (see, e.g., \cite{Mena2018} for the
differential Lyapunov and Riccati matrix equations). In fact,
it can be shown that projecting the Turing solution on a low-rank
manifold, especially during the transient unstable time dynamics,
can induce the selection of specific Fourier modes in the final
pattern (see the discussion above). This topic will be object of
future investigations.
The right plot of Figure~\ref{fig:SchnackSolCost} 
reports the computational costs of all methods. We recall that by
varying  $N_x$, Sylvester matrix equations in \eqref{reacdiff_matrix} of increasing 
size are solved by the reduced spectral approach. As in 
Test (a), rEuler and rSBDF have almost the same cost and
are cheaper than the other methods and ADI is less expensive than 
rExp. {For the largest spatial dimensions rEuler becomes the most effective,
costs-wise.}
{We stress} that for the larger values of 
$N_x$ this test could not be performed by classical vector-oriented version of 
the same schemes due to the high computational load.

\begin{figure}[tb]
\centering
\includegraphics[scale=0.32]{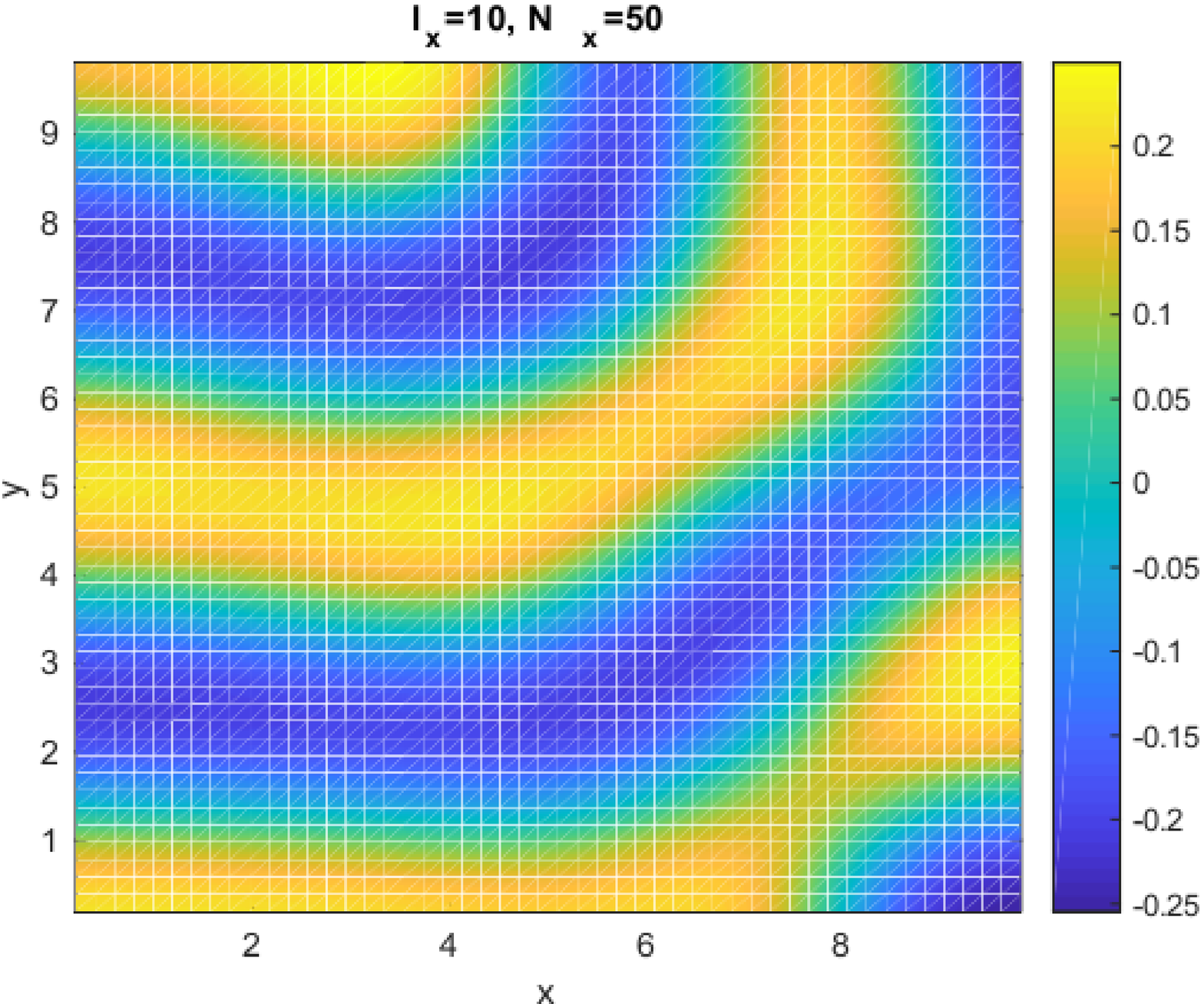}
\includegraphics[scale=0.33]{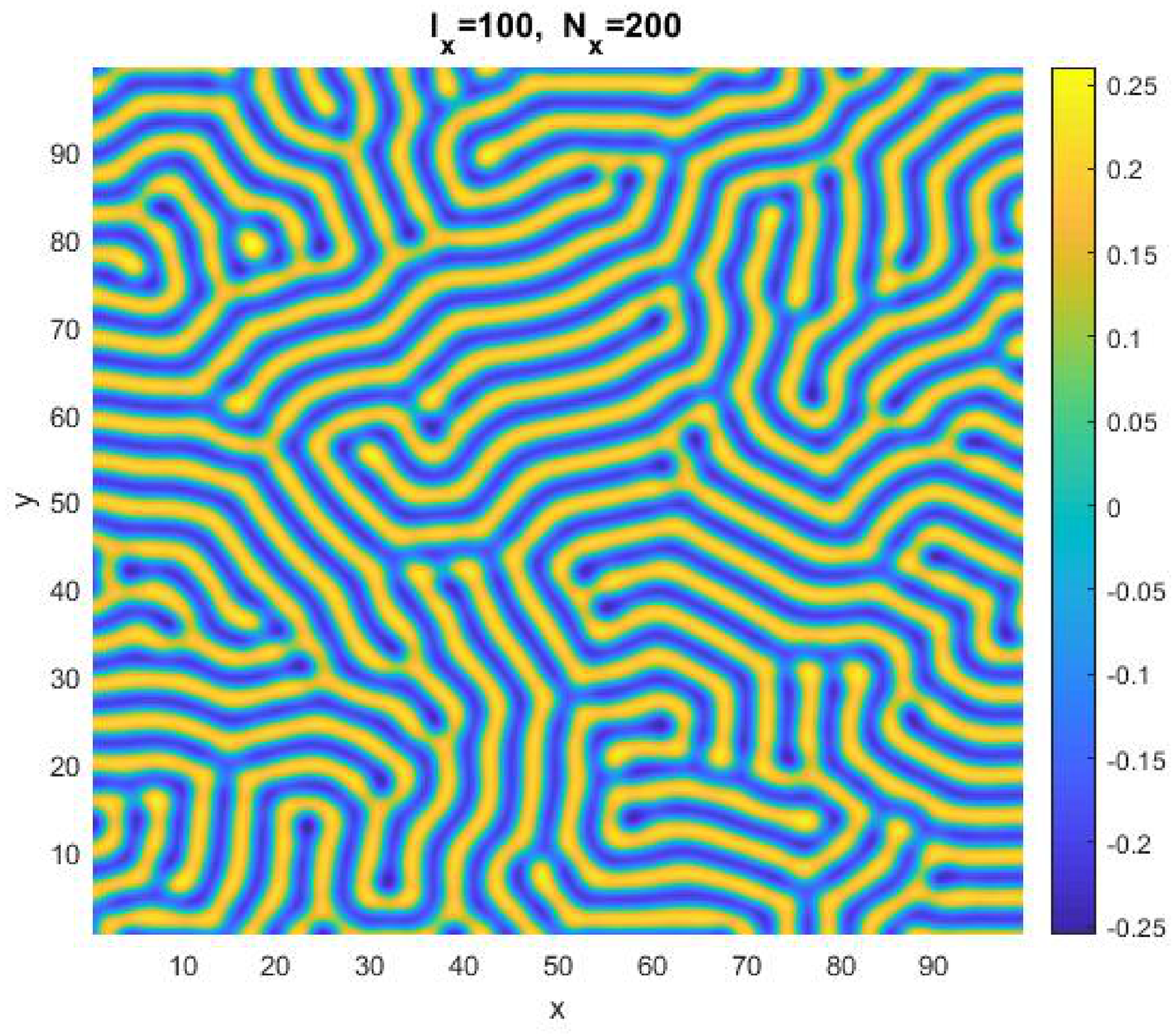}
\includegraphics[scale=0.33]{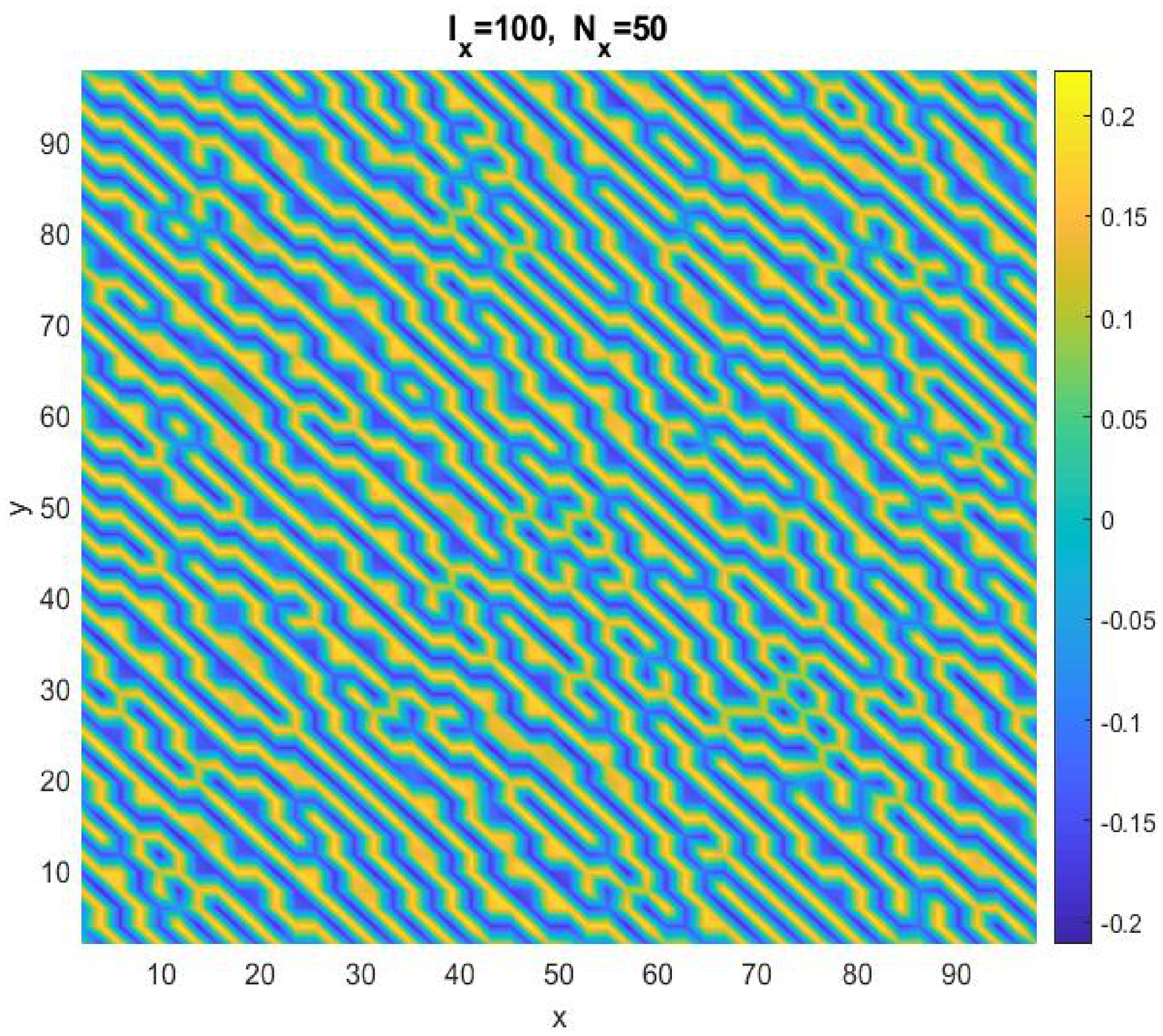}
\caption{Labyrinth Turing pattern of DIB Model. Top left: $\Omega=[0,10]\times [0,10]$ and
$N_x=50 (h_x=0.2)$.
Top right: $\Omega=[0,100]\times [0,100]$ and $N_x=200$ $(h_x=0.5)$. Bottom: $\Omega=[0,100]\times [0,100]$ and $N_x=50$ $(h_x=2)$.}
\label{fig:DIB_lab}
\end{figure}
\subsection{DIB model}
{In this section we report on} the importance of 
the matrix-oriented approach to carefully approximate the spatial 
structure of Turing patterns on fine meshgrids and large domains at feasible computational costs.
\begin{figure}[htb]
\centering
\includegraphics[scale=0.39]{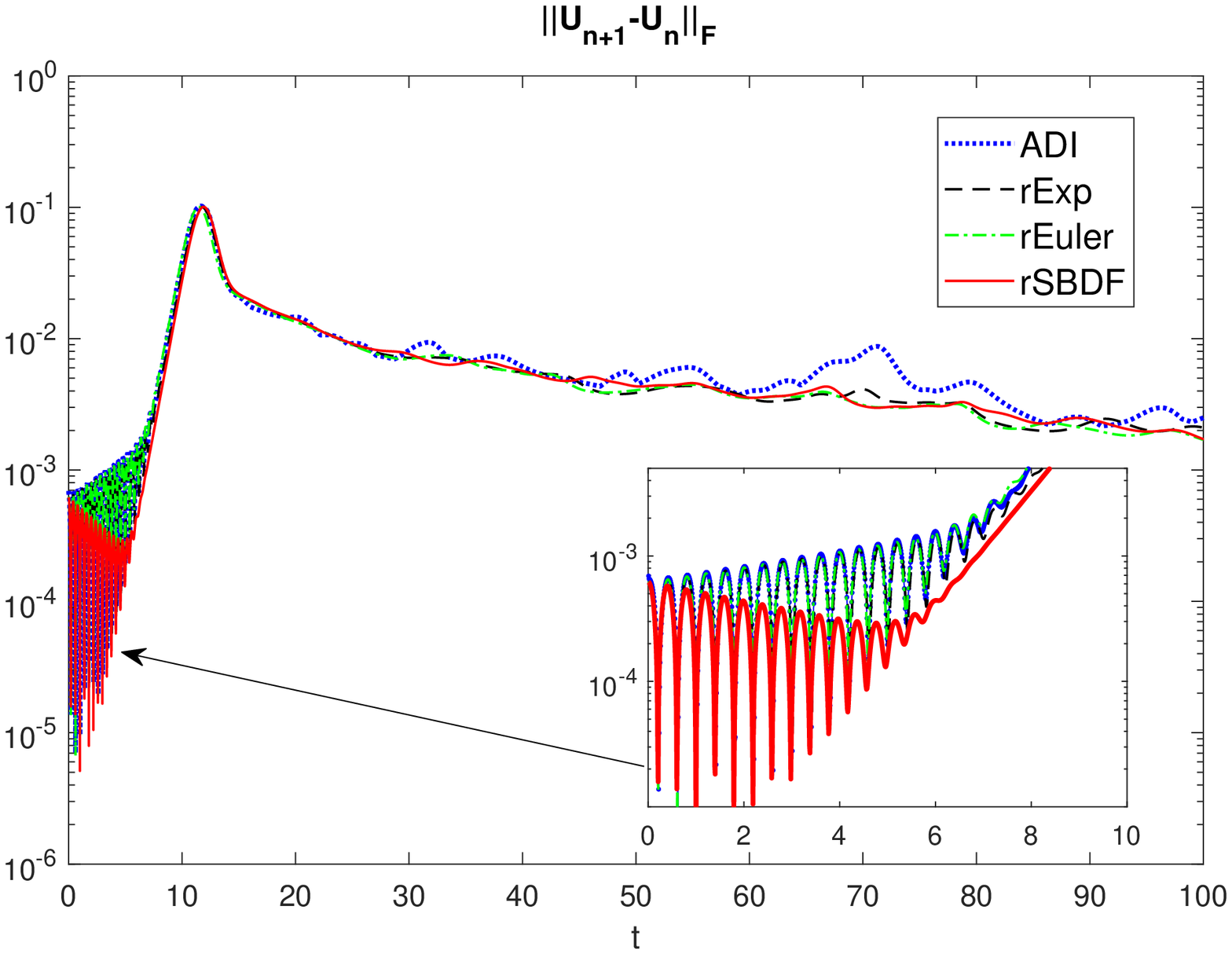}
\includegraphics[scale=0.39]{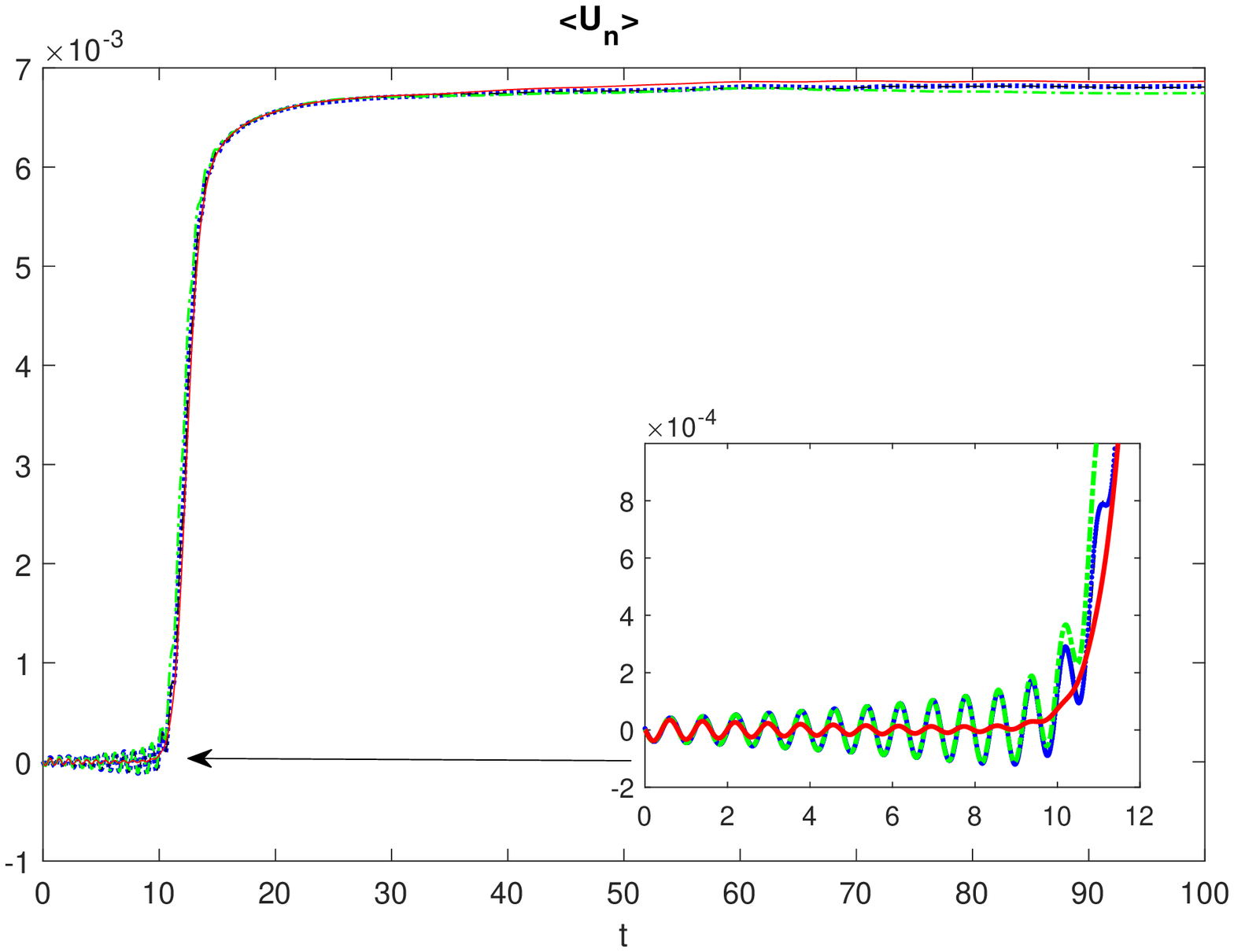}
\caption{DIB model: time dynamics of the increment $\delta_n = \|U_{n+1}-U_n \|_F$ (left plot) and of the mean value $<U_n >$ (right plot) for all methods in the case $\Omega=[0,100]\times [0,100]$, $T_f=100$, $h_t=10^{-2}$, $N_x=200$.}
\label{DIB}
\end{figure}
Towards this aim we consider the RD-PDE model studied in 
\cite{EJAM2015,IPSE2018} describing an electrodeposition process for 
metal growth where the kinetics in \eqref{reacdiff} are given by
\begin{equation}
\begin{split}
&f_1(u,v)=\rho \left( A_1(1-v)u-A_2 \,u^3-B(v-\alpha) \right),\\
&f_2(u,v)=\rho \left(C (1+k_2 u)(1-v)[1-\gamma(1-v)]-D v (1+k_3 u)(1+\gamma v))\right).
\label{CineticheDIB}
\end{split}\end{equation}
Here $u(x,y,t)$ represents the morphology of the metal deposit, while $v(x,y,t)$ 
monitors its surface percentual chemical composition. 
The nonlinear source terms account for generation and loss of relevant material 
during the process. In \cite{Lacit2017} this model has been 
proposed to study pattern formation during the charge-discharge process of batteries. 
In the same article, it has also been proved that for a given parameter choice 
of the RD-PDE model there exists an \emph{intrinsic
 pattern type} that can only emerge if an 
\emph{effective domain size} of integration is considered, and this is
 given by ${\cal A} = \rho |\Omega|$, 
where $|\Omega|$ = area($\Omega$). Hence, if the scaling factor 
in \eqref{CineticheDIB} is $\rho=1$, a large domain $\Omega$ must be 
chosen to ``see'' the Turing pattern, and the grid fineness sufficiently high to
capture the pattern details. For this reason, the number of 
meshpoints $N_x,N_y$, that is the size of the Sylvester equations \eqref{reacdiff_matrix}, must 
be sufficiently large. {Figure~\ref{fig:DIB_lab} reports three typical situations:
the upper left plot refers to a too small domain to be able to identify the 
morphological class,
which is instead clearly visible in the upper right plot, determined with a much larger domain and
a fine grid. In the third setting (large domain but a too coarse grid, $N_x=50$) shown in the lower plot
of Figure~\ref{fig:DIB_lab}, the numerical approximation is unable to  clearly detect
the labyrinth pattern in its full granularity.}
These solutions have been obtained by solving \eqref{reacdiff}-\eqref{CineticheDIB} on a
square domain $\Omega=[0, \ell_x] \times [0, \ell_x]$ and with the following parameter choice for which a labyrinth
pattern is expected (\cite{Lacit2017,EJAM2015,IPSE2018}):
\begin{eqnarray*}
d_1=1, d_2=20, \rho=1, A_1=10, A_2=30, k_2=2.5,k_3=1.5, 
\alpha=0.5, \gamma=0.2, D=2.4545, B=66, C=3.
\end{eqnarray*}

\noindent We have applied {ADI} and the matrix methods rEuler, rExp, rSBDF until $T_f=100$, with $h_t=10^{-2}$. 
In Figure~\ref{DIB}} we show the dynamics of the
increment $\delta_n$ and of the mean value $<U_n>$ for the simulations corresponding to the full labyrinth in Figure~\ref{fig:DIB_lab}(upper right). For the chosen (large) $h_t$ the methods exhibit different reactivity and stabilizing properties. {The rSBDF method seems to display} the best performance.
\begin{table}[bht]
\centering
\begin{tabular}[]{|c|r|r|r|}
\hline
\textbf{Methods} & $\ell_x=10, N_x=50$ & $\ell_x=100, N_x=200$ & $\ell_x=100, N_x=300$\\
\hline
\hline
IMEX Euler (vector form) & 3.5742 &  90.2079  &  234.7923 \\
\hline
rEuler &2.1592   & 34.3328 & 79.6804 \\
\hline
rSBDF &  3.3318  &  51.1874 & 118.4675 \\
\hline
rExp & 4.6030  & 50.6238 & 127.8123 \\
\hline
ADI & 3.2780  & 44.3830 & 95.0798 \\
\hline
\end{tabular}
\caption{DIB model. Computational times in seconds for all methods to obtain the two top patterns in Figure~\ref{fig:DIB_lab}.
The cost for the case $N_x=300, \ell_x=100$ is also reported.
\label{tab:cputimes}}
\end{table}

{Table~\ref{tab:cputimes} reports the computational times of all
numerical methods for obtaining the patterns in Figure~\ref{fig:DIB_lab}, 
that is for the two cases (i) $N_x=50, \ell_x =10,$ (left) and 
(ii) $N_x=200, \ell_x=100$ (right), including that of the IMEX Euler vector formulation 
(implemented similarly to that used in Table~\ref{tab:Heat}). The cost for a larger spatial meshgrid $N_x=300$ on $\Omega=[0,100] \times [0,100]$ is also reported.
For $N_x=50$, that is when the pattern is not well identified, all methods 
display similar computational performances. 
In the other cases, the vector form significantly suffers from dealing with much larger dimensional
data, with respect to the matrix-oriented schemes.
rEuler exhibits the best computational times for all dimensions, requiring about
one third of the other methods' time for the large values of $N_x$.
The other matrix-based schemes are almost equivalent, with ADI being slightly less expensive. 
These preliminary experiments seem to indicate that the matrix formulation is a competitive
methodology  for the numerical solution 
of the RD-PDE systems when a fine spatial grid is necessary to capture the morphological features of the pattern. }

\section{Conclusions}
{By exploiting the Kronecker structure of the diffusion matrix, we have shown that
the classical semi-discretization in space of reaction-diffusion PDEs on regular domains can 
be seen as a system of {\it matrix} ODE equations (see \eqref{eqn:matrixform}). 
ODE solvers in time such as IMEX Euler, 2SBDF schemes and the Exponential Euler method,
 can thus be implemented in matrix form, requiring a sequence
of small matrix problems (Sylvester equations, matrix exponentials) to be solved at each
time step. Due to the modest size of these matrices, the computational cost per iteration 
can be made lower than that of the corresponding vector approaches, by working
in the (reduced) spectral space. 
To avoid the high computational load of the vector-IMEX methods, in the literature
this challenge has often been faced by the using ADI approach;
our comparisons show that the new \emph{reduced} schemes can be a valid alternative to ADI. 
In particular, rEuler exhibits the best performance. To the best of our knowledge, the exponential method rExp
is new in this field of application.
The improvement obtained by working with matrix ODEs allows us to capture the key features of Turing pattern solutions, and this
has been shown by using two typical benchmarks, the Schnackenberg and DIB models. 
We plan to deepen our understanding of matrix-oriented formulation by further exploring  higher
order methods so as to improve the accuracy of the methods while maintaining efficiency.
Finally, we speculate that the matrix-approach can be particularly helpful, for instance, in the context 
of parameter identification problems (see e.g. \cite{IPSE2018}) to significantly reduce 
the computational load for the corresponding constrained minimization problem, 
which requires the approximations of PDE models at each optimization step.
\section*{Acknowledgements}
The work of the third author was partially supported by the Indam-GNCS 2017 Project
``Metodi numerici avanzati per equazioni e funzioni di matrici con struttura'' and by
the Project AlmaIdea, Universit\`a di Bologna. All authors are
members of the GNCS-INdAM activity group. 

\section*{Appendix: Matlab code}
{
This Appendix reports the possible implementation in Matlab (\cite{matlab}) of the reduced Euler Algorithm
and of the reduced Exponential Euler Algorithm for a RD-PDE system of two equations.
The displayed codes compute and employ explicit inverse matrices. This procedure turned out
to be more effective than solving the corresponding systems on the fly at each iteration.
We also note that the reported {\bf rEuler} code explicitly computes the eigenvalue decomposition
of $A_1, A_2$ and $B$, which could be avoided: in {\bf rExp} we showed how to derive the
eigendecompositions of all matrices from that of $T$.

{\footnotesize
\begin{verbatim}
%%% INPUT:
% F1, F2 = kinetics of the RD-model (inputs to be edited as needed)
% T  = Approximation of the 1D laplacian, order p, of size Nx-1
% Lx = length of spatial domain: [0,Lx]x[0,Lx]
% Nx = number of points for the spatial discretization in each direction
% Tf = final time of integration
% ht = time step
% p = order of approximation of the second derivative
% par = parameters structure for the RD-model
% w = correction coefficient (only for rExp)

%%% OUPUT:
% U, V = solution of the models
%%%%%%%%%%%%%%%%%%%%%%%%%%%%%%%%%%%%%%%%%%%%%%%% rEuler
function [U,V]=rEuler(F1, F2, T, Lx, Nx, Tf, ht, p, ue, ve, d, par) 
hs = Lx/(Nx); 
I = eye(Nx-1);
Nt = Tf/ht;
A1 = I-ht*T;             A2=I-d*ht*T;
[Qa1,Ra1] = eig(A1);     [Qa2,Ra2] = eig(A2); 
I_Qa1 = inv(Qa1);        I_Qa2 = inv(Qa2);
B = -ht*T';            [Qb,Rb]=eig(B);       
                       I_Qb = inv(Qb);
dA1 = diag(Ra1); dA2 = diag(Ra2); dB = diag(Rb);
Inveig_u = 1./(dA1*ones(1,Nx-1) + ones(Nx-1,1)*dB');
Inveig_v = 1./(dA2*ones(1,Nx-1) + ones(Nx-1,1)*d*dB');
rng('default');
Rs = rand(Nx-1); %perturbation matrix for the IC
U = ue+Rs*1e-05;
V = ve+Rs*1e-05;
for k=1:Nt 
     C_u = U+ht*F1(U, V, par);
     C_v = V+ht*F2(U, V, par);
     Cc_u = I_Qa1*C_u*Qb;
     Cc_v = I_Qa2*C_v*Qb;
     Xx = Cc_u.*Inveig_u;
     Yy = Cc_v.*Inveig_v;
     U = Qa1*Xx*I_Qb;
     V = Qa2*Yy*I_Qb;
end

%%%%%%%%%%%%%%%%%%%%%%%%%%%%%%%%%%%%%%%%%%%%%%%% rExp
function [U,V]=rExp(F1, F2, T, Lx, Nx, Tf, ht, p, ue, ve, d, par, w) 
hs = Lx/(Nx); 
I = eye(Nx-1);
Nt = Tf/ht;
rng('default');
Rs = rand(Nx-1); %perturbation matrix for the IC
U1 = ue+Rs*1e-05; 
V1 = ve+Rs*1e-05;

[QT,ET]=eig(T);
A1plus = T-w*I;
Q1plus = QT;       I_Q1plus = inv(QT);      E1plus = ET-w*I;
Bplus = T'; 
QBplus = inv(QT'); I_QBplus = QT';          EBplus = conj(ET); 
A2plus = d*T - w*I;
Q2plus = QT;       I_Q2plus = inv(Q2plus);  E2plus = d*ET-w*I;
EAB1 = exp(ht*diag(E1plus))*exp(ht*diag(EBplus).');
EAB2 = exp(ht*diag(E2plus))*exp(ht*d*diag(EBplus).');
LL1 = 1./(ht*diag(E1plus)*ones(1,Nx-1)+ht*ones(Nx-1,1)*diag(EBplus).');
LL2 = 1./(ht*diag(E2plus)*ones(1,Nx-1)+ht*d*ones(Nx-1,1)*diag(EBplus).');
U1eig = Q1plus\U1*QBplus;
V1eig = Q2plus\V1*QBplus;
for k=1:Nt 
    G1 = F1(U1,V1, par) + w*U1;
    ProjG1 = I_Q1plus*G1*QBplus;
    ErhsU = EAB1.*ProjG1-ProjG1;
    U = EAB1.*(U1eig)+ ht*(LL1.*ErhsU);  
    U1eig = U;
    U = Q1plus*U*I_QBplus;
    
    G2 = F2(U1,V1, par)+w*V1; 
    ProjG2 = I_Q2plus*G2*QBplus;
    ErhsV = EAB2.*ProjG2-ProjG2;
    V = EAB2.*V1eig + ht*(LL2.*ErhsV);
    V1eig = V;
    V = Q2plus*V*I_QBplus;
    U1 = U; V1 = V;
end
\end{verbatim}
}
}


\end{document}